\documentclass[11pt]{article}

\usepackage{amssymb,latexsym}
\usepackage{amsmath,amscd}
\usepackage{citesort}
\usepackage{theorem}
\usepackage{eucal}
\usepackage{bbm}

\title{Completely positive inner products\\ 
  and\\
  strong Morita equivalence}

\author{\textbf{Henrique Bursztyn}\thanks{E-mail:
    henrique@math.toronto.edu}\addtocounter{footnote}{5}
  \\[0.1cm]
  Department of Mathematics\\
  University of Toronto\\
  100 St. George Street\\
  Toronto, Ontario, M5S 3G3\\
  Canada\\[0.5cm]
  \textbf{Stefan Waldmann}\thanks{E-mail:
    Stefan.Waldmann@physik.uni-freiburg.de}
  \\[0.1cm]
  Fakult{\"a}t f{\"u}r Mathematik und Physik\\
  Albert-Ludwigs-Universit{\"a}t Freiburg\\
  Physikalisches Institut\\
  Hermann Herder Stra{\ss}e 3\\
  D 79104 Freiburg\\
  Germany}

\date{October 2004}

\renewcommand{\mathbb}[1]{\mathbbm{#1}} 

\newcommand{\textdef}[1] {\textbf{#1}}

\newcommand{\im}         {{\mathrm i}}
\newcommand{\eu}         {{\mathrm e}}

\newcommand{\cc} [1]     {\overline{{#1}}}
\newcommand{\id}         {{\mathsf{id}}}
\newcommand{\tr}         {{\mathrm{tr}}}

\newcommand{\ad}         {{\mathrm{ad}}}
\newcommand{\Aut}        {{\mathrm{Aut}}}

\newcommand{\Hom}        {{\mathrm{Hom}}}
\newcommand{\Mor}        {{\mathrm{Mor}}}
\newcommand{\End}        {{\mathrm{End}}}

\newcommand{\st} [1]     {\scriptscriptstyle{{#1}}}
\newcommand{\Ad}         {{\mathrm{Ad}}}

\newcommand{\SP} [1]     {{\left\langle{{#1}}\right\rangle}}
\newcommand{\Unit}       {\mathbb{1}}

\newcommand{\Exp}        {{\mathrm{Exp}}}
\newcommand{\Ln}         {{\mathrm{Ln}}}
\newcommand{\ring}[1]    {{\mathsf{{#1}}}}
\newcommand{\HdR}        {\mathrm{H}_{\st{\mathrm{dR}}}}

\newcommand{\Ped}[1]    {{\mathcal{P}(\mathcal{#1})}}

\newcommand{\Bimod}[3] {{\sideset{_{\scriptscriptstyle{#1}}}{_{\scriptscriptstyle{#3}}}{\operatorname{#2}}}}

\newcommand{\BEA}    {\Bimod{\mathcal{B}}{\mathcal{E}}{\mathcal{A}}}   
\newcommand{\FEA}    {\Bimod{\mathfrak{F}(\mathcal{E})}{\mathcal{E}}{\mathcal{A}}}
\newcommand{\AccEB}  {\Bimod{\mathcal{A}}{\cc{\mathcal{E}}}{\mathcal{B}}} 
\newcommand{\AEpB}   {\Bimod{\mathcal{A}}{\mathcal{E}^\prime}{\mathcal{B}}}

\newcommand{\FEpB}   {\Bimod{\mathfrak{F}(\mathcal{E}^\prime)}{\mathcal{E}}{\mathcal{B}}}
\newcommand{\BccEpF} {\Bimod{\mathcal{B}}{\cc{\mathcal{E}^\prime}}{\mathfrak{F}(\mathcal{E}^\prime)}}

\newcommand{\AAA}   {\Bimod{\mathcal{A}}{\mathcal{A}}{\mathcal{A}}}   
\newcommand{\BBB}   {\Bimod{\mathcal{B}}{\mathcal{B}}{\mathcal{B}}}   
\newcommand{\PBEPA} {\Bimod{\mathcal{P}_{\mathcal{B}}}{\mathcal{E}}{\mathcal{P}_{\mathcal{A}}}}

\newcommand{\DDD}   {\Bimod{\mathcal{D}}{\mathcal{D}}{\mathcal{D}}}
\newcommand{\CFB}   {\Bimod{\mathcal{C}}{\mathcal{F}}{\mathcal{B}}}
\newcommand{\BpEA}  {\Bimod{\mathcal{B}}{\mathcal{E}^\prime}{\mathcal{A}}}
\newcommand{\CpFB}  {\Bimod{\mathcal{C}}{\mathcal{F}^\prime}{\mathcal{B}}}
\newcommand{\AHD}   {\Bimod{\mathcal{A}}{\mathcal{H}}{\mathcal{D}}}

\newcommand{\FB}    {\Bimod{}{\mathcal{F}}{\mathcal{B}}}
\newcommand{\GC}    {\Bimod{}{\mathcal{G}}{\mathcal{C}}}
\newcommand{\DGDp}  {\Bimod{\mathcal{D}}{\mathcal{G}}{\mathcal{D}^\prime}}

\newcommand{\cstBEA} {\Bimod{\mathcal{B}}{\widehat{\mathcal{E}}}{\mathcal{A}}}
\newcommand{\cstCFB} {\Bimod{\mathcal{C}}{\widehat{\mathcal{F}}}{\mathcal{B}}}
\newcommand{\cstBFA} {\Bimod{\mathcal{B}}{\widehat{\mathcal{F}}}{\mathcal{A}}}

\newcommand{\tensor}[1][{}]{\mathbin{\otimes_{\scriptscriptstyle{\mathcal{{#1}}}}}}
\newcommand{\tensM}[1][{}] {\mathbin{\widehat{\otimes}_{\scriptscriptstyle{\mathcal{{#1}}}}}}
\newcommand{\tensB}[1][{}] {\mathbin{\widetilde{\otimes}_{\scriptscriptstyle{\mathcal{{#1}}}}}}
\newcommand{\tensMF}[1][{}] {\mathbin{\widehat{\otimes}_{\scriptscriptstyle{{{#1}}}}}}

\newcommand{\repC}[1][{}] {\sideset{^*}{_{\mathcal{#1}}}{\operatorname{\textrm{-}{\mathsf{rep}}}}}
\newcommand{\RepC}[1][{}] {\sideset{^*}{_{\mathcal{#1}}}{\operatorname{\textrm{-}{\mathsf{Rep}}}}}
\newcommand{\rep}[1][{}]  {\sideset{^*}{_{\mathcal{#1}}}{\operatorname{\textrm{-}\mathsf{rep}}}}

\newcommand{\modrep}[1][{}]  {\sideset{^*}{_{\mathcal{#1}}}{\operatorname{\textrm{-}\mathsf{mod}}}}
\newcommand{\Modrep}[1][{}]  {\sideset{^*}{_{\mathcal{#1}}}{\operatorname{\textrm{-}\mathsf{Mod}}}}

\newcommand{\Pic}  {\mathsf{Pic}}
\newcommand{\StrPic} {\mathsf{Pic}^{\mathrm{str}}}
\newcommand{\starPic} {\mathsf{Pic}^*}

\newcommand{\catstar} {{^*}\mbox{-}\mathsf{Alg}}
\newcommand{\catstarpos} {{^*}\mbox{-}\mathsf{Alg}^+}
\newcommand{\Cstar} {\mathsf{C}^*}

\newcommand{\IP}[4]{{\,}_{\scriptscriptstyle{#2}\!\!}\left\langle{{#1}}\right\rangle^{\scriptscriptstyle{#3}}_{\scriptscriptstyle{#4}}}

\newcommand{\SPEA}[1] {\IP{{#1}}{}{\mathcal{E}}{\mathcal{A}}}
\newcommand{\BSPE}[1] {\IP{{#1}}{\mathcal{B}}{\mathcal{E}}{}}
\newcommand{\SPFB}[1] {\IP{{#1}}{}{\mathcal{F}}{\mathcal{B}}}

\newcommand{\SPEcstA}[1] {\IP{{#1}}{}{\widehat{\mathcal{E}}}{\mathcal{A}}}
\newcommand{\BSPEcst}[1] {\IP{{#1}}{\mathcal{B}}{\widehat{\mathcal{E}}}{}}

\newcommand{\SPFEA}[1] {\IP{{#1}}{}{\mathcal{F}\otimes\mathcal{E}}{\mathcal{A}}}
\newcommand{\BSPFE}[1] {\IP{{#1}}{\mathcal{B}}{\mathcal{F}\otimes\mathcal{E}}{}}

\newcommand{\SPA}[1] {\IP{{#1}}{}{}{\mathcal{A}}}
\newcommand{\ASP}[1] {\IP{{#1}}{\mathcal{A}}{}{}}

\newtheorem{lemma} {Lemma} [section]
\newtheorem{proposition}[lemma]{Proposition}

\newtheorem{theorem}[lemma]{Theorem}
\newtheorem{corollary}[lemma]{Corollary}
\newtheorem{definition}[lemma]{Definition}

\theoremstyle{definition}
\newtheorem{example}[lemma]{Example}
\newtheorem{remark}[lemma]{Remark}

\newenvironment{proof}{{\sc Proof:}}{{\hspace*{\fill} $\square$\\}}

\numberwithin{equation}{section}

\begin{document}

\maketitle

\begin{abstract}
    We develop a general framework for the study of strong Morita
    equivalence in which $C^*$-algebras and hermitian star products on
    Poisson manifolds are treated in equal footing. We compare strong
    and ring-theoretic Morita equivalences in terms of their Picard
    groupoids for a certain class of unital $^*$-algebras encompassing
    both examples.  Within this class, we show that both notions of
    Morita equivalence induce the same equivalence relation but
    generally define different Picard groups. For star products, this
    difference is expressed geometrically in cohomological terms.
\end{abstract}


%
%

\tableofcontents

%
%
\section{Introduction}
\label{sec:intro}

This paper investigates several similarities between two types of
algebras with involution: hermitian star products on Poisson manifolds
and $C^*$-algebras.  Their connection is suggested by their common
role as ``quantum algebras'' in mathematical physics, despite the fact
that the former is a purely algebraic notion, whereas the latter has
important analytical features.  Building on
\cite{bursztyn.waldmann:2001a,bursztyn.waldmann:2001b}, we develop in
this paper a framework for their unified study, focusing on Morita
theory; in particular, the properties shared by $C^*$-algebras and
star products allow us to develop a general theory of \textit{strong}
Morita equivalence in which they are treated in equal footing.

Our set-up is as follows.  We consider $^*$-algebras over rings of the
form $\ring{C}=\ring{R}(\im)$, where $\ring{R}$ is an ordered ring and
$\im^2=-1$. The main examples of $\ring{R}$ that we will have in mind
are $\mathbb{R}$, with its natural ordering, and
$\mathbb{R}[[\lambda]]$, with ordering induced by ``asymptotic
positivity'', i.e., $a = \sum_{r=0}^\infty a_r \lambda^r > 0$ if and
only if $a_{r_0}>0$, where $a_{r_0}$ is the first nonzero coefficient
of $a$. This general framework encompasses complex $^*$-algebras, such
as $C^*$-algebras, as well as $^*$-algebras over the ring of formal
power series $\mathbb{C}[[\lambda]]$, such as hermitian star products.
We remark that the case of $^*$-algebras over $\mathbb{C}$ has
been extensively studied, see e.g. \cite{schmuedgen:1990a}, and
\cite{waldmann:2003c:pre} for a comparison of notions of
positivity.

In our general framework, we define a purely algebraic notion of
strong Morita equivalence.  The key ingredient in this definition
is the notion of \textit{completely positive} inner products,
which we use to refine Ara's $^*$-Morita equivalence
\cite{ara:1999a}.  One of our main results is that completely
positive inner products behave well under the internal and
external tensor products, and, as a consequence, strong Morita
equivalence defines an equivalence relation within the class of
non-degenerate and idempotent $^*$-algebras. This class of
algebras includes both star products and $C^*$-algebras as
examples.  We prove that important constructions in the theory of
$C^*$-algebras, such as Rieffel's induction of representations
\cite{rieffel:1974a}, carry over to this purely algebraic setting,
recovering and improving many of our previous results
\cite{bursztyn.waldmann:2001b,bursztyn.waldmann:2001a}.

In the ordinary setting of unital rings, Morita equivalence
coincides with the notion of isomorphism in the category whose
objects are unital rings and morphisms are isomorphism classes of
bimodules, composed via tensor product.  The invertible arrows in
this category form the \textit{Picard groupoid} $\Pic$
\cite{benabou:1967a}, which is a ``large'' groupoid (in the sense
that its collection of objects is not a set) encoding the
essential aspects of Morita theory: the orbit of a ring in $\Pic$
is its Morita equivalence class, whereas the isotropy groups in
$\Pic$ are the usual Picard groups of rings.
Analogously, we show that our purely algebraic notion of strong
Morita equivalence coincides with the notion of isomorphism in a
category whose objects are non-degenerate and idempotent
$^*$-algebras over a fixed ring $\ring{C}$; morphisms and their
compositions are given by more elaborate bimodules and tensor
products, and  the invertible arrows in this category form the
\textit{strong Picard groupoid} $\StrPic$.
When restricted to $C^*$-algebras, we show that $\StrPic$ defines
an equivalence relation which turns out to coincide with Rieffel's
(analytical) notion of strong Morita equivalence
\cite{rieffel:1974b}, and its isotropy groups are the Picard
groups of $C^*$-algebras as in \cite{brown.green.rieffel:1977a};
these results are proven along the lines of
\cite{ara:1999b,bursztyn.waldmann:2001b}.

In the last part of the paper, we compare strong and
ring-theoretic Morita equivalences for unital $^*$-algebras over
$\ring{C}$ by analyzing the canonical groupoid morphism
\begin{equation}
    \label{eq:morpintro}
    \StrPic \longrightarrow \Pic.
\end{equation}
We prove that, for a suitable class of unital $^*$-algebras,
including both unital $C^*$-algebras and hermitian star products,
$\StrPic$ and $\Pic$ have the same orbits, i.e., the two notions
of Morita equivalence define the \textit{same} equivalence
relation. This is a simultaneous extension of Beer's result
\cite{beer:1982a}, in the context of $C^*$-algebras, and
\cite[Thm.~2]{bursztyn.waldmann:2002a}, for deformation
quantization. Despite the coincidence of orbits, we show that, for
both unital $C^*$-algebras and hermitian star products, the
isotropy groups of $\Pic$ and $\StrPic$ are generally different.
We note that the obstructions for \eqref{eq:morpintro} being an
equivalence can be described in a unified way for both classes of
$^*$-algebras, due to common properties of their automorphism
groups. A key ingredient for this discussion in the context of
formal deformation quantization is the fact that hermitian star
products are always (completely) \textit{positive} deformations,
in the sense that positive measures on the manifold can be
deformed into positive linear functionals of the star product, see
\cite{bursztyn.waldmann:2000a,bursztyn.waldmann:2004}.

The paper is organized as follows: In Section~\ref{sec:pos} we
recall the basic definitions and properties of $^*$-algebras over
ordered rings and pre-Hilbert spaces.
Section~\ref{sec:InnerProducts} is devoted to completely positive
inner products, a central notion throughout the paper. In
Section~\ref{sec:TensorProducts} we define various categories of
representations and prove that internal and external tensor
products of completely positive inner products are again
completely positive. In Section~\ref{sec:StrMorEq} we define
strong Morita equivalence, prove that it is an equivalence
relation within the class of non-degenerate and idempotent
$^*$-algebras and show that strong Morita equivalence implies the
equivalence of the categories of representations introduced in
Section~\ref{sec:TensorProducts}.  In Section~\ref{sec:Picard} we
define the strong Picard groupoid and relate our algebraic
definition to the $C^*$-algebraic Picard groupoid, proving their
equivalence. In Section~\ref{sec:strongversusalegbraic} we study
the map \eqref{eq:morpintro} for a suitable class of unital
$^*$-algebras.
Finally, in Section~\ref{sec:applications},
we consider hermitian deformations, and, in particular, hermitian
star products.

\noindent \textbf{Acknowledgments:} We would like to thank Martin
Bordemann for many valuable discussions. We also thank S. Jansen, N.
Landsman, I. Moerdijk, R.  Nest, and A. Weinstein for useful comments
and remarks. H.B.  thanks DAAD for financial support and Freiburg University and IPAM-UCLA for their
hospitality while part of this work was being done.

\noindent
\textbf{Conventions}: Throughout this paper $\ring{C}$ will denote a
ring of the form $\ring{R}(\im)$, where $\ring{R}$ is an ordered ring
and $\im^2 = -1$.  Unless otherwise stated, algebras and modules will
always be over a fixed ring $\ring{C}$.  For a manifold $M$,
$C^\infty(M)$ denote its algebra of \textit{complex-valued} smooth
functions.

%
%

\section{$^*$-Algebras, positivity and pre-Hilbert spaces}
\label{sec:pos}

A \textdef{$^*$-algebra} over $\ring{C}$ is a $\ring{C}$-algebra
equipped with an anti-linear involutive anti-automorphism.  If
$\mathcal{A}$ is a $^*$-algebra over $\ring{C}$, then there are
natural notions of positivity induced by the ordering structure on
$\ring{R}$: A \textdef{positive linear functional} is a
$\ring{C}$-linear map $\omega: \mathcal{A} \longrightarrow
\ring{C}$ satisfying $\omega(a^*a) \ge 0$ for all $a \in
\mathcal{A}$, and an algebra element $a \in \mathcal{A}$ is called
\textdef{positive} if $\omega(a) \ge 0$ for all positive linear
functionals $\omega$ of $\mathcal{A}$. Elements of the form
\begin{equation}
    \label{eq:pos}
    r_1 a_1^*a_1 + \cdots + r_n a_n^*a_n,
\end{equation}
$r_i \in \ring{R}^+$, $a_i \in \mathcal{A}$ are clearly positive.
The set of positive algebra elements is denoted by
$\mathcal{A}^+$, see \cite[Sec.~2]{bursztyn.waldmann:2001a} for
details.  These definitions recover the standard notions of
positivity when $\mathcal{A}$ is a $C^*$-algebra; for
$\mathcal{A}=C^\infty(M)$, positive linear functionals coincide
with positive Borel measures on $M$ with compact support, and
positive elements are positive functions
\cite[App.~B]{bursztyn.waldmann:2001a}.

A linear map $\phi:\mathcal{A} \longrightarrow \mathcal{B}$, where
$\mathcal{A}$ and $\mathcal{B}$ are $^*$-algebras over $\ring{C}$, is
called \textdef{positive} if $\phi(\mathcal{A}^+) \subseteq
\mathcal{B}^+$, and \textdef{completely positive} if the canonical
extensions $\phi: M_n(\mathcal{A}) \longrightarrow M_n(\mathcal{B})$
are positive for all $n \in \mathbb{N}$.
\begin{example}
    \label{ex:cpmaps}
    Let us consider the maps $\tr:M_n(\mathcal{A}) \longrightarrow
    \mathcal{A}$ and $\tau: M_n(\mathcal{A}) \longrightarrow
    \mathcal{A}$ defined by
    \begin{equation}
        \label{eq:TraceAndTau}
        \tr(A) = \sum_{i=1}^n A_{ii},
        \quad\textrm{and}\quad
        \tau(A) = \sum_{i,j=1}^n A_{ij},
    \end{equation}
    where $A = (A_{ij}) \in M_n(\mathcal{A})$. A direct computation shows
    that both maps are positive. Replacing $\mathcal{A}$ by
    $M_N(\mathcal{A})$ and using the identification
    $M_n(M_N(\mathcal{A})) \cong M_{Nn}(\mathcal{A})$, it immediately
    follows that $\tr$ and $\tau$ are completely positive maps.
\end{example}

A \textdef{pre-Hilbert space} $\mathcal{H}$ over $\ring{C}$ is a
$\ring{C}$-module with a $\ring{C}$-valued sesquilinear inner product
satisfying
\begin{equation}
    \label{eq:InnerProdHilbertSpace}
    \SP{\phi,\psi} = \cc{\SP{\psi,\phi}}
    \quad
    \textrm{and}
    \quad
    \SP{\phi,\phi} > 0
    \;
    \textrm{for}
    \;
    \phi \ne 0,
\end{equation}
see \cite{bursztyn.waldmann:2001a}. We use the convention that
$\SP{\cdot,\cdot}$ is linear in the second argument. These are
direct analogues of complex pre-Hilbert spaces. A
\textdef{$^*$-representation} of a $^*$-algebra $\mathcal{A}$ on a
pre-Hilbert space $\mathcal{H}$ is a $^*$-homomorphism from
$\mathcal{A}$ into the adjointable endomorphisms
$\mathfrak{B}(\mathcal{H})$ of $\mathcal{H}$, see
\cite{bursztyn.waldmann:2001a,bursztyn.waldmann:2001b}; the main
examples are the usual representations of $C^*$-algebras on
Hilbert spaces and the formal representations of star products,
see e.g. \cite{bordemann.waldmann:1998a,waldmann:2002a}.

%
%

\section{Completely positive inner products}
\label{sec:InnerProducts}

%
%

\subsection{Inner products and complete positivity}
\label{subsec:InnProdCompPos}

Let $\mathcal{A}$ be a $^*$-algebra over $\ring{C}$ and consider a
right $\mathcal{A}$-module $\mathcal{E}$. Thoughout this paper,
$\mathcal{A}$-modules are always assumed to have a compatible
$\ring{C}$-module structure.
\begin{remark}
    \label{rem:convent}
    When $\mathcal{A}$ is unital, we adopt the convention that $x
    \cdot \Unit = x$ for $x \in \mathcal{E}$; morphisms between unital
    algebras are assumed to be unital.
\end{remark}

An \textdef{$\mathcal{A}$-valued inner product} on $\mathcal{E}$
is a $\ring{C}$-sesquilinear (linear in the second argument) map
$\SP{\cdot,\cdot}: \mathcal{E} \times \mathcal{E} \longrightarrow
\mathcal{A}$ so that, for all $x, y \in \mathcal{E}$ and $a \in
\mathcal{A}$,
\begin{equation}
    \label{eq:InnerProduct}
    \SP{x,y} = \SP{y,x}^*
    \quad
    \textrm{and}
    \quad
    \SP{x, y \cdot a} = \SP{x,y} a.
\end{equation}
The definition of an $\mathcal{A}$-valued inner product on a left
$\mathcal{A}$-module is analogous, but we require linearity on the
first argument. We call an inner product $\SP{\cdot,\cdot}$
\textdef{non-degenerate} if $\SP{x,y} = 0$ for all $y$ implies that $x
= 0$, and \textdef{strongly non-degenerate} if the map $\mathcal{E}
\longrightarrow \Hom_{\mathcal{A}}(\mathcal{E},\mathcal{A})$, $x
\mapsto \SP{x,\cdot}$ is a bijection.  Two inner products
$\SP{\cdot,\cdot}_1$ and $\SP{\cdot,\cdot}_2$ on $\mathcal{E}$ are
called \textdef{isometric} if there exists a module automorphism $U$
with $\SP{Ux,Uy}_1 = \SP{x,y}_2$.

An endomorphism $T \in \End_{\mathcal{A}}(\mathcal{E})$ is
\textdef{adjointable} with respect to $\SP{\cdot,\cdot}$ if there
exists $T^* \in \End_{\mathcal{A}}(\mathcal{E})$ (called an
\textdef{adjoint} of $T$) such that
\begin{equation}
    \label{eq:Adjoint}
    \SP{x, Ty} = \SP{T^*x,y}
\end{equation}
for all $x,y \in \mathcal{E}$. The algebra of adjointable
endomorphisms is denoted by $\mathfrak{B}_{\mathcal{A}}(\mathcal{E})$,
or simply $\mathfrak{B}(\mathcal{E})$. If $\SP{\cdot,\cdot}$ is
non-degenerate, then adjoints are unique and
$\mathfrak{B}_{\mathcal{A}}(\mathcal{E})$ becomes a $^*$-algebra over
$\ring{C}$. One defines the $\ring{C}$-module
$\mathfrak{B}_{\mathcal{A}}(\mathcal{E}, \mathcal{F})$ of adjointable
homomorphisms $\mathcal{E} \longrightarrow \mathcal{F}$ analogously.

We now use the positivity notions in $\mathcal{A}$: An inner
product $\SP{\cdot,\cdot}$ on $\mathcal{E}$ is \textdef{positive}
if $\SP{x,x} \in \mathcal{A}^+$ for all $x \in \mathcal{E}$, and
\textdef{positive
  definite} if $0 \ne \SP{x,x} \in \mathcal{A}^+$ for $x \ne 0$.
\begin{definition}
    \label{definition:CompletelyPositive}
    Consider $\mathcal{E}^n$ as a right $M_n(\mathcal{A})$-module,
    and let $\SP{\cdot,\cdot}^{(n)}$ be the $M_n(\mathcal{A})$-valued
    inner product on $\mathcal{E}^n$ defined by
    \begin{equation}
        \label{eq:SpnExtension}
        \SP{x,y}^{(n)}_{i j} = \SP{x_i,y_j},
    \end{equation}
    where $x = (x_1, \ldots, x_n)$ and $y = (y_1, \ldots, y_n) \in
    \mathcal{E}^n$. We say that $\SP{\cdot,\cdot}$ is
    \textdef{completely positive} if $\SP{\cdot,\cdot}^{(n)}$ is
    positive for all $n$.
\end{definition}
\begin{remark}
    \label{remark:DirectSums}
    Although the direct sum of non-degenerate (resp. positive,
    completely positive) inner products is non-degenerate (resp.
    positive, completely positive), this may not hold for positive
    definiteness: Consider $\mathcal{A} = \mathbb{Z}_2$ as
    $^*$-algebra over $\mathbb{Z}(\im)$, see
    \cite[Sec.~2]{bursztyn.waldmann:2001b}; then the canonical inner
    product on $\mathcal{A}$ is positive definite but on
    $\mathcal{A}^2$ the vector $(\Unit,\Unit)$ satisfies
    $\SP{(\Unit,\Unit), (\Unit,\Unit)} = \Unit + \Unit = 0$.
\end{remark}

The following  observation provides a way to detect algebras
$\mathcal{A}$ for which positive $\mathcal{A}$-valued inner products
on arbitrary $\mathcal{A}$-modules are automatically completely
positive.
\begin{proposition}
    \label{proposition:cpmatrix}
    Let $\mathcal{A}$ be a $^*$-algebra satisfying the following
    property: for any $n \in \mathbb{N}$, if $(A_{ij}) \in
    M_n(\mathcal{A})$ satisfies $\sum_{ij} a_i^*A_{ij}a_j \in
    \mathcal{A}^+$ for all $(a_1,\ldots,a_n) \in \mathcal{A}^n$, then
    $A \in M_n(\mathcal{A})^+$.  Then any positive
    $\mathcal{A}$-valued inner product on an $\mathcal{A}$-module is
    automatically completely positive.
\end{proposition}
\begin{proof}
    Let $\mathcal{E}$ be an $\mathcal{A}$-module with positive inner
    product $\SP{\cdot,\cdot}$, and let $x_1,\ldots, x_n \in
    \mathcal{E}$. For $a_1,\ldots,a_n \in \mathcal{A}$, the matrix
    $A=(\SP{x_i,x_j})$ satisfies
    \[
    \sum\nolimits_{ij}a_i^*\SP{x_i,x_j}a_j
    = \sum\nolimits_{ij}\SP{x_i \cdot a_i, x_j \cdot a_j}
    = \SP{
      \sum\nolimits_i x_i \cdot a_i,
      \sum\nolimits_j x_j \cdot a_j}
    \in \mathcal{A}^+.
    \]
    So the matrix $(\SP{x_i,x_j})$ is positive and $\SP{\cdot,\cdot}$
    is completely positive.
\end{proof}
The converse also holds, e.g., if $\mathcal{A}$ is unital.

Note that, although a positive definite inner product is always
non-degenerate, a positive inner product which is non-degenerate
may fail to be positive definite. This is due to the fact that the
\textdef{degeneracy space} of an $\mathcal{A}$-module
$\mathcal{E}$ with inner product $\SP{\cdot,\cdot}$, defined by
\begin{equation}
    \label{eq:DegeneracySpace}
    \mathcal{E}^\bot
    = \{x \in \mathcal{E} \; | \; \SP{x,\cdot} = 0\},
\end{equation}
might be \emph{strictly} contained in the space
\begin{equation}
    \label{eq:normzerospace}
    \{x \in \mathcal{E} \; | \; \SP{x,x} = 0\}.
\end{equation}
\begin{example}
    \label{example:NonAdmissible}
    Let $\mathcal{A}=\bigwedge^\bullet(\ring{C}^n)$ be the Grassmann
    algebra over $\ring{C}^n$, with $^*$-involution defined by $e_i^*
    = e_i$, where $e_1, \ldots, e_n$ is the canonical basis for
    $\ring{C}^n$.  Regard $\mathcal{A}$ as a right module over itself,
    equipped with inner product $\SP{x,y}=x^*\wedge y$. Then $\SP{e_i,e_i} = 0$. However, $\mathcal{A}^\bot = \{0\}$,
    since $\SP{1, x} = x$.
\end{example}

Any $\mathcal{A}$-valued inner product on $\mathcal{E}$ induces a
non-degenerate one on the quotient $\mathcal{E} \big/
\mathcal{E}^\bot$. Moreover, (completely) positive inner products
induce (completely) positive inner products.  In case
$\mathcal{E}^\bot = \{x \in \mathcal{E} \; | \; \SP{x,x}=0\}$, the
quotient inner product is \textit{positive definite}.

$^*$-Algebras possessing a ``large'' amount of positive linear
functionals, such as $C^*$-algebras and formal hermitian deformation
quantizations \cite{bursztyn.waldmann:2001b,bursztyn.waldmann:2001a},
are such that \eqref{eq:DegeneracySpace} and \eqref{eq:normzerospace}
coincide.
\begin{example}
    \label{example:SuffManyPosAdmissible}
    Let $\mathcal{A}$ be a $^*$-algebra over $\ring{C}$ with the
    property that, for any non-zero hermitian element $a \in
    \mathcal{A}$, there exists a positive linear functional $\omega$
    with $\omega(a)\neq 0$. Under the additional assumption that $2
    \in \ring{R}$ is invertible, any $\mathcal{A}$-module
    $\mathcal{E}$ with $\mathcal{A}$-valued inner product is such that
    \eqref{eq:DegeneracySpace} and \eqref{eq:normzerospace} coincide.
    The proof follows from the arguments in
    \cite[Sect.~5]{bursztyn.waldmann:2001a}.
\end{example}

%
%
\subsection{Examples of completely positive inner products}
\label{subsec:ExamplesCP}

Inner products on complex pre-Hilbert spaces are always completely
positive. This result extends in two directions: on one hand, one can
replace $\mathbb{C}$ by arbitrary rings $\ring{C}$; on the other hand,
$\mathbb{C}$ can be replaced by more general $C^*$-algebras.
\begin{example}
    \label{example:ComPosPreHilbert}
    \emph{(Pre-Hilbert spaces over $\ring{C}$)}
    \\
    If $\mathcal{A} = \ring{C}$, then
    \cite[Prop.~A.4]{bursztyn.waldmann:2001a} shows that the condition
    in Proposition~\ref{proposition:cpmatrix} is satisfied. So a
    positive $\ring{C}$-valued inner product on any $\ring{C}$-module
    $\mathcal{H}$ is completely positive.  This is the case, in
    particular, for inner products on pre-Hilbert spaces over
    $\ring{C}$ (which are non-degenerate).
\end{example}

\begin{example}
    \label{example:ComPosCstarAlg}
    \emph{(Pre-Hilbert $C^*$-modules)}
    \\
    Let $\mathcal{A}$ be a $C^*$-algebra over $\ring{C} = \mathbb{C}$.
    Then the condition in Proposition~\ref{proposition:cpmatrix}
    holds, see e.g. \cite[Lem.~2.28]{raeburn.williams:1998a}.  So a
    positive $\mathcal{A}$-valued inner product $\SP{\cdot,\cdot}$ on
    any $\mathcal{A}$-module $\mathcal{E}$ is completely positive, see
    also \cite[Lem.~2.65]{raeburn.williams:1998a}.  When
    $\SP{\cdot,\cdot}$ is positive definite,
    $(\mathcal{E},\SP{\cdot,\cdot})$ is called a \textdef{pre-Hilbert
      $C^*$-module} over $\mathcal{A}$.
\end{example}

Example~\ref{example:ComPosPreHilbert} uses the quotients fields
of $\ring{R}$ and $\ring{C}$, whereas
Example~\ref{example:ComPosCstarAlg} uses the functional calculus
of $C^*$-algebras, so neither immediately extend to inner products
with values in arbitrary $^*$-algebras. Nevertheless, one can
still show the complete positivity of particular inner products.

\begin{example}
    \label{example:FreeCanonical}
    \emph{(Free modules)}
    \\
    Consider $\mathcal{A}^N$ as a right $\mathcal{A}$-module with
    respect to right multiplication, equipped with the canonical inner
    product
    \begin{equation}
        \label{eq:CanonicalInnerProduct}
        \SP{x,y} = \sum_{i=1}^N x_i^*y_i,
    \end{equation}
    where $x=(x_1,\ldots,x_N), y=(y_1,\ldots,y_N) \in \mathcal{A}^N$.
    This inner product is completely positive since, for $x^{(1)},
    \ldots, x^{(n)} \in \mathcal{A}^N$, the matrix $X =
    \left(\SP{x^{(\alpha)}, x^{(\beta)}}\right) \in M_n(\mathcal{A})$
    can be written as
    \begin{equation}
        \label{eq:XsumXi}
        X =  \sum_{i} X_i^* X_i,
        \quad
        \textrm{where}
        \quad
        X_i =
        \left(
            \begin{matrix}
                x_i^{(1)} & \ldots & x_i^{(n)} \\
                0 & \ldots & 0 \\
                \vdots & & \vdots \\
                0 & \ldots & 0
            \end{matrix}
        \right).
    \end{equation}
    Note, however, that the inner product
    \eqref{eq:CanonicalInnerProduct} need not be positive definite as
    there may exist elements $a \in \mathcal{A}$ with $a^*a = 0$.  If
    $\mathcal{A}$ is unital then \eqref{eq:CanonicalInnerProduct} is
    strongly non-degenerate; in the non-unital case it may be
    degenerate.
\end{example}
\begin{remark}
    \label{rem:Pinner}
    Let $\mathcal{E}$ be an $\mathcal{A}$-module with inner product
    $\SP{\cdot,\cdot}$ which can be written as
    \begin{equation}
        \label{eq:Pinner}
        \SP{x,y}=\sum_{i=1}^m P_i(x)^*P_i(y), \;\; \mbox{ for } x,y \in
        \mathcal{E},
    \end{equation}
    where $P_i:\mathcal{E} \longrightarrow
    \mathcal{A}$ are $\mathcal{A}$-linear maps. By replacing
    $x_i^{(\alpha)}$ with $P_i(x^{(\alpha)})$ in
    Example~\ref{example:FreeCanonical}, one immediately sees that
    \eqref{eq:Pinner} is completely positive.
\end{remark}

A direct computation shows that completely positive inner products
restrict to completely positive inner products on submodules.
\begin{example}
    \label{ex:Submodules}
    \emph{(Hermitian projective modules)}
    \\
    The restriction of the canonical inner product
    \eqref{eq:CanonicalInnerProduct} to any submodule of
    $\mathcal{A}^n$ is completely positive.  In particular, hermitian
    projective modules, i.e., modules of the form $\mathcal{E} =
    P\mathcal{A}^n$, where $P \in M_n(\mathcal{A})$, $P= P^2 = P^*$,
    have an induced completely positive inner product (this also follows from
    Remark~\ref{rem:Pinner}). If $\mathcal{A}$ is
    \emph{unital}, then this inner product is strongly non-degenerate.
\end{example}

The following simple observation concerns uniqueness.
\begin{lemma}
    \label{lem:unique}
    Let $\mathcal{E}$ be an $\mathcal{A}$-module equipped with a
    strongly non-degenerate $\mathcal{A}$-valued inner product
    $\SP{\cdot,\cdot}$. Let $\SP{\cdot,\cdot}'$ be another inner
    product on $\mathcal{E}$. Then there exists a unique hermitian
    element $H \in \mathfrak{B}(\mathcal{E})$ such that
    \begin{equation}
        \label{eq:OtherInnProdh}
        \SP{x,y}' = \SP{x, Hy},
    \end{equation}
    and $\SP{\cdot,\cdot}'$ is isometric to $\SP{\cdot,\cdot}$ if
    there exists an invertible $U \in \mathfrak{B}(\mathcal{E})$ with
    $H = U^*U$.
\end{lemma}

\begin{example}
    \label{example:HermVecBundle}
    \emph{(Hermitian vector bundles)}
    \\
    Let $\mathcal{A}=C^\infty(M)$ be the algebra of smooth
    complex-valued functions on a manifold $M$. As a result of
    Serre-Swan's theorem \cite{swan:1962a}, hermitian projective
    modules $P\mathcal{A}^N$
    correspond to (sections of) vector bundles over $M$ (since in this
    case idempotents are always equivalent to projections, see
    Section~\ref{subsec:StructBimod}), and $\mathcal{A}$-valued inner
    products correspond to hermitian fiber metrics.

    As noticed in Example~\ref{ex:Submodules}, there is a strongly
    non-degenerate inner product $\SP{\cdot,\cdot}$ on
    $P\mathcal{A}^N$.  For any other inner product
    $\SP{\cdot,\cdot}'$, there exists a unique hermitian element $H
    \in PM_N(\mathcal{A})P$ such that $\SP{x,y}' = \SP{x, Hy}$.  Since
    any positive invertible element $H \in PM_N(C^\infty(M))P$ can be
    written as $H = U^*U$ for an invertible $U\in PM_N(C^\infty(M))P$,
    it follows from Lemma \ref{lem:unique} that there is only one
    fiber metric on a vector bundle over $M$ up to isometric
    isomorphism. We will generalize this example in
    Section~\ref{subsec:StructBimod}.
\end{example}
\begin{example}
    \label{example:NontrivialInnProd}
    \emph{(Nontrivial inner products)}
    \\
    Even if the algebra $\mathcal{A}$ is a field, one can
    have nontrivial inner products.  For example, consider $\ring{R} =
    \mathbb{Q}$ and $\ring{C} = \mathbb{Q}(\im)$. Then $3 \in
    \ring{C}$ is a positive invertible element but there is no $z \in
    \ring{C}$ with $\cc{z}z = 3$ (write $z = a + \im b$ with $a =
    r/n$, $b=s/n$ with $r,s,n\in \mathbb{N}$, then take the equation
    $3n^2 = r^2 + s^2$ modulo $4$). Hence $\SP{z,w}' = 3\cc{z}w$ is
    completely positive and strongly non-degenerate but not isometric
    to the canonical inner product $\SP{z,w} = \cc{z}w$.
\end{example}

%
%

\section{Representations and tensor products}
\label{sec:TensorProducts}

%
%

\subsection{Categories of $^*$-representations}
\label{subsec:CatRepTheories}

We now discuss the algebraic analogues of Hilbert $C^*$-modules, see
e.g. \cite{lance:1995a}.  Let $\mathcal{D}$ be a $^*$-algebra over
$\ring{C}$.
\begin{definition}
    \label{definition:PreHilbertMod}
    A (right) \textdef{inner-product $\mathcal{D}$-module} is a pair
    $(\mathcal{H},\SP{\cdot,\cdot})$, where $\mathcal{H}$ is a (right)
    $\mathcal{D}$-module and $\SP{\cdot,\cdot}$ is a non-degenerate
    $\mathcal{D}$-valued inner product. If $\SP{\cdot,\cdot}$ is
    completely positive, we call $(\mathcal{H},\SP{\cdot,\cdot})$ a
    \textdef{pre-Hilbert $\mathcal{D}$-module}.
\end{definition}
Whenever there is no risk of confusion, we will denote an
inner-product module (or pre-Hilbert module) simply by $\mathcal{H}$.

We now consider $^*$-representations of $^*$-algebras on inner-product
modules, extending the discussion in
\cite{bursztyn.waldmann:2001a,bursztyn.waldmann:2001b,bordemann.waldmann:1998a}.
Let $\mathcal{A}$ be a $^*$-algebra over $\ring{C}$, and let
$\mathcal{H}$ be an inner-product $\mathcal{D}$-module.
\begin{definition}
    \label{definition:StarReps}
    A \textdef{$^*$-representation} of $\mathcal{A}$ on $\mathcal{H}$
    is a $^*$-homomorphism $\pi: \mathcal{A} \longrightarrow
    \mathfrak{B}_{\mathcal{D}}(\mathcal{H})$.
\end{definition}
An \textdef{intertwiner} between two $^*$-representations
$(\mathcal{H}, \pi)$ and $(\mathcal{K}, \varrho)$ is an isometry $T
\in \mathfrak{B}_{\mathcal{D}}(\mathcal{H}, \mathcal{K})$ such that,
for all $a \in \mathcal{A}$,
\begin{equation}
    \label{eq:Intertwiner}
    T \pi(a) = \varrho(a) T.
\end{equation}

We denote by $\modrep[D](\mathcal{A})$ the category whose objects
are $^*$-representations of $\mathcal{A}$ on inner-product modules
over $\mathcal{D}$ and morphisms are intertwiners. The subcategory
whose objects are $^*$-representations on pre-Hilbert modules is
denoted by $\repC[D](\mathcal{A})$. Since both categories contain
trivial representations of $\mathcal{A}$, we will consider the
following further refinement:
A $^*$-representation $(\mathcal{H},\pi)$ is \textdef{strongly
  non-degenerate} if
\begin{equation}
    \label{eq:nondeg}
    \pi(\mathcal{A})\mathcal{H} = \mathcal{H},
\end{equation}
(by Remark \ref{rem:convent}, this is always the case if $\mathcal{A}$ is
unital). The category of strongly non-degenerate $^*$-representations
of $\mathcal{A}$ on inner-product (resp.  pre-Hilbert) $\mathcal{D}$-modules
is denoted by $\Modrep[D](\mathcal{A})$ (resp.
$\RepC[D](\mathcal{A})$).
\begin{definition}
    \label{definition:InnerProductBimodule}
    An inner-product $\mathcal{D}$-module $\mathcal{H}$ together with
    a strongly non-degenerate $^*$-re\-pre\-sen\-tation of $\mathcal{A}$
    will be called an
    $(\mathcal{A},\mathcal{D})$-\textdef{inner-product bimodule}; it is
    an $(\mathcal{A},\mathcal{D})$-\textdef{pre-Hilbert
      bimodule} if $\mathcal{H}$ is a pre-Hilbert module.
\end{definition}
These are algebraic analogues of Hilbert bimodules as e.g. in
\cite[Def.~3.2]{landsman:2001b}. This terminology differs from the one in
\cite{ara:1999b}.

An \textdef{isomorphism} of inner-product (bi)modules (or pre-Hilbert
(bi)modules) is just a (bi)module homomorphism preserving inner
products.

More generally, suppose $\AHD$ is a bimodule equipped with an arbitrary
$\mathcal{D}$-valued inner product $\SP{\cdot,\cdot}$. We say that
$\SP{\cdot,\cdot}$ is \textdef{compatible} with the
$\mathcal{A}$-action if
\begin{equation}
    \label{eq:compatible}
    \SP{a\cdot x, y} = \SP{x,a^*\cdot y},
\end{equation}
for all $a \in \mathcal{A}$ and $x,y \in \mathcal{H}$. Clearly, any
$^*$-representation of $\mathcal{A}$ on an inner-product module
$\mathcal{H}$ over $\mathcal{D}$ makes $\mathcal{H}$ into a bimodule
for which $\SP{\cdot,\cdot}$ and the $\mathcal{A}$-action are
compatible. Unless otherwise stated, inner products on bimodules are
assumed to be compatible with the actions.

%
%

\subsection{Tensor products and Rieffel induction of representations}
\label{subsec:RieffelInduction}

Let $\mathcal{A}$ and $\mathcal{B}$ be $^*$-algebras over $\ring{C}$.
Let $\FB$ be a right $\mathcal{B}$-module equipped with a
$\mathcal{B}$-valued inner product $\SPFB{\cdot,\cdot}$, and let
$\BEA$ be a bimodule equipped with an $\mathcal{A}$-valued inner
product $\SPEA{\cdot,\cdot}$ compatible with the $\mathcal{B}$-action.
Following Rieffel \cite{rieffel:1974a,rieffel:1974b}, there is a
well-defined $\mathcal{A}$-valued inner product $\SPFEA{\cdot,\cdot}$
on the tensor product $\FB \tensor[B] \BEA$ completely determined by
\begin{equation}
    \label{eq:InternalTensorProduct}
    \SPFEA{y_1 \otimes_{\st{\mathcal{B}}} x_1, y_2 \otimes_{\st{\mathcal{B}}} x_2}
    = \SPEA{x_1, \SPFB{y_1, y_2} \cdot x_2}
\end{equation}
for $x_1, x_2 \in \BEA$ and $y_1, y_2 \in \FB$ (we extend it to
arbitrary elements using $\ring{C}$-sesquilinearity).  An analogous
construction works for left modules carrying inner products.

If $(\FB \tensor[B] \BEA)^\bot$ is the degeneracy space associated
with $\SPFEA{\cdot,\cdot}$, then the quotient
\[
(\FB \tensor[B] \BEA) \big/ (\FB \tensor[B] \BEA)^\bot
\]
acquires an induced inner product, also denoted by
$\SPFEA{\cdot,\cdot}$, which is non-degenerate, see
Section~\ref{subsec:InnProdCompPos}. Thus the pair
\begin{equation}
    \label{eq:Internal}
    \FB \tensM[B] \BEA := \left(
        (\FB \tensor[B] \BEA) \big/ (\FB \tensor[B] \BEA)^\bot,
        \SPFEA{\cdot,\cdot}
    \right)
\end{equation}
is an inner-product $\mathcal{A}$-module called the \textdef{internal
  tensor product} of $(\FB, \SPFB{\cdot,\cdot})$ and
$(\BEA,\SPEA{\cdot,\cdot})$.  As we will see, in many examples the
degeneracy space of \eqref{eq:InternalTensorProduct} is already
trivial.
\begin{lemma}
    \label{lemma:Degenerate}
    If $\mathcal{C}$ is a $^*$-algebra and $\CFB$ is a bimodule so
    that $\SPFB{\cdot,\cdot}$ is compatible with the
    $\mathcal{C}$-action, then $\FB \tensM[B] \BEA$ carries a
    canonical left $\mathcal{C}$-action, compatible with
    $\SPFEA{\cdot,\cdot}$.
\end{lemma}
The proof of this lemma is a direct computation. It is also simple to
check that internal tensor products have associativity properties
similar to those of ordinary (algebraic) tensor products: Let $\GC$ be
a $\mathcal{C}$-module with $\mathcal{C}$-valued inner product, and
let $\CFB$ (resp. $\BEA$) be a bimodule with $\mathcal{B}$-valued
(resp.  $\mathcal{A}$-valued) inner product compatible with the
$\mathcal{C}$-action (resp. $\mathcal{B}$-action).
\begin{lemma}
    \label{lem:assoc}
    There is a natural isomorphism
    \begin{equation}
        \label{eq:NaturalAssocIso}
        (\GC \tensM[C] \CFB) \tensM[B] \BEA
        \cong
        \GC \tensM[C] (\CFB \tensM[B] \BEA)
    \end{equation}
    induced from the usual associativity of algebraic tensor products.
\end{lemma}

Internal tensor products also behave well with respect to maps.
\begin{lemma}
    \label{lemma:MorphInternalTensorMorph}
    Let $\CFB, \CpFB$ be equipped with compatible $\mathcal{B}$-valued
    inner products, and let $\BEA$, $\BpEA$ be equipped with
    compatible $\mathcal{A}$-valued inner products.  Let $S \in
    \mathfrak{B}(\CFB, \CpFB)$ and $T \in \mathfrak{B}(\BEA, \BpEA)$
    be adjointable bimodule morphisms.  Then their algebraic tensor
    product $S \tensor[B] T$ induces a well-defined adjointable
    bimodule morphism $S \tensM[B] T: \CFB \tensM[B] \BEA
    \longrightarrow \CpFB \tensM[B] \BpEA$ with adjoint given by $S^*
    \tensM[B] T^*$. If $S$ and $T$ are isometric then $S \tensM[B] T$
    is isometric as well.
\end{lemma}

Hence, for a fixed triple of $^*$-algebras $\mathcal{A}$,
$\mathcal{B}$ and $\mathcal{C}$, we obtain a \emph{functor}
\begin{equation}
    \label{eq:tensMBFunctor}
    \tensM[B]: \modrep[B](\mathcal{C}) \times \modrep[A](\mathcal{B})
    \longrightarrow \modrep[A](\mathcal{C}).
\end{equation}
In the case of unital algebras, one can replace $\modrep$ by
$\Modrep$ in \eqref{eq:tensMBFunctor}.

A central question is whether one can restrict the functor
$\tensM[B]$ to representations on \emph{pre-Hilbert modules}, or
whether the tensor product \eqref{eq:InternalTensorProduct} of two
\emph{positive} inner products remains positive. This is the case,
for example, in the realm of $C^*$-algebras, but the proof uses
the functional calculus, see e.g.
\cite[Prop.~2.64]{raeburn.williams:1998a}. Fortunately, a purely
algebraic result can be obtained if one requires the inner
products to be \emph{completely positive}.
\begin{theorem}
    \label{theorem:CPTensorCP}
    If the inner products $\SPFB{\cdot,\cdot}$ on $\FB$ and
    $\SPEA{\cdot,\cdot}$ on $\BEA$ are completely positive, then the
    inner product $\SPFEA{\cdot,\cdot}$ on $\FB \tensor[B] \BEA$
    defined by \eqref{eq:InternalTensorProduct} is also completely
    positive.
\end{theorem}
\begin{proof}
    Let $\Phi^{(1)}, \ldots \Phi^{(n)} \in \FB \tensor[B]
    \BEA$. We must show that the matrix
    $$
    A =
    \left(\SPFEA{\Phi^{(\alpha)}, \Phi^{(\beta)}}\right)
    $$
    is a positive element in $M_n(\mathcal{A})$.  Without loss of
    generality, we can write $\Phi^{(\alpha)} = \sum_{i=1}^N
    y_i^{(\alpha)} \tensor[B] x_i^{(\alpha)}$, where $N$
    is the same for all $\alpha$.  Let us consider the map
    \begin{equation}
        \label{eq:YAPM}
        f: M_{n}(M_N(\mathcal{B}))
        \longrightarrow
        M_{n}(M_N(\mathcal{A})),\quad
        (B_{ij}^{\alpha\beta})
        \mapsto
        \left(
            \SPEA{x_i^{(\alpha)},
              B_{ij}^{\alpha\beta} \cdot x^{(\beta)}_j}
        \right),
    \end{equation}
    $1\leq i,j \leq N$, $1\leq \alpha,\beta \leq n$.  We claim that
    $f$ is a positive map. Indeed, as a consequence of the definition
    of positive maps in Sect.~\ref{sec:pos}, it suffices to show that
    $f(B^*B)$ is positive for any $B \in M_n(M_N(\mathcal{B}))$. A
    direct computation shows that, for $B = (B_{ij}^{\alpha\beta})$,
    $$
    f(B^*B)= \sum_{k=1}^N\sum_{\gamma=1}^n C_k^\gamma
    \quad
    \textrm{with}
    \quad
    \left(C_k^\gamma\right)^{\alpha\beta}_{ij}
    =
    \SPEA{B_{ki}^{\gamma\alpha} x_i^{(\alpha)},
      B_{kj}^{\gamma\beta}x_j^{(\beta)}},
    $$
    which is positive since $\SPEA{\cdot,\cdot}$ is completely
    positive.  Since the matrix
    $$
    \left(\SPFB{y^{(\alpha)}_i,
          y^{(\beta)}_j}\right) \in M_{nN}(\mathcal{B})
    $$
    is positive, for $\SPFB{\cdot,\cdot}$ in $\mathcal{F}$ is completely
    positive, it follows that the matrix
    $$
    \left(\SPEA{x_i^{(\alpha)},
          \SPFB{y^{(\alpha)}_i, y^{(\beta)}_j} \cdot
          x^{(\beta)}_j}\right)
    $$
    is a positive element in $M_{nN}(\mathcal{A})$.  Since
    summation over $i,j$ defines a positive map $\tau:
    M_{nN}(\mathcal{A}) \longrightarrow M_n(\mathcal{A})$,  see
    Example~\ref{ex:cpmaps}, the matrix
    \begin{equation}
        \label{eq:YAPMatrix}
        \sum_{i,j=1}^N
        \left(
            \SPEA{x_i^{(\alpha)},
              \SPFB{y^{(\alpha)}_i,
                y^{(\beta)}_j} \cdot x^{(\beta)}_j}
        \right)
        = \left(\SPFEA{\Phi^{(\alpha)}, \Phi^{(\beta)}}\right)
        = A
    \end{equation}
    is positive. This concludes the proof.
\end{proof}

As pointed out in Section~\ref{subsec:InnProdCompPos}, if
$\SPFEA{\cdot,\cdot}$ is completely positive, then so is the
induced inner product on $\FB \tensM[B] \BEA$.
\begin{corollary}
    \label{corollary:tenscompos}
    If $\FB$ and $\BEA$ have completely positive inner products, then
    $\FB \tensM[B] \BEA$ is a pre-Hilbert module.
\end{corollary}

It follows that the functor $\tensM[B]$ in
\eqref{eq:tensMBFunctor} restricts to a functor
\begin{equation}
    \label{eq:tensMBFunctor2}
    \tensM[B]: \rep[B](\mathcal{C}) \times \rep[A](\mathcal{B})
    \longrightarrow \rep[A](\mathcal{C}),
\end{equation}
and, from \eqref{eq:tensMBFunctor2}, we obtain two functors by
fixing each one of the two arguments.
\begin{example}
    \label{example:Rieffel}
    \emph{(Rieffel induction)}
    \\
    Let $\mathcal{A}$, $\mathcal{B}$ and $\mathcal{D}$ be
    $^*$-algebras, and fix a
    $(\mathcal{B},\mathcal{A})$-bimodule $\BEA \in
    \repC[A](\mathcal{B})$.  We then have a functor
    \begin{equation}
        \label{eq:RieffelInduction}
        \mathsf{R}_{\mathcal{E}} = \BEA \tensM[A] \cdot :
        \repC[D](\mathcal{A})
        \longrightarrow \repC[D](\mathcal{B});
    \end{equation}
    on objects, $\mathsf{R}_{\mathcal{E}}(\AHD) = \BEA \tensM[A]
    \AHD$, and, on morphisms,
    $\mathsf{R}_{\mathcal{E}}(T)= \id \tensM[A] T$, for $T\in
    \mathfrak{B}(\mathcal{H},\mathcal{H}')$.

    This functor is called \textdef{Rieffel induction} and relates the
    representation theories of $\mathcal{A}$ and $\mathcal{B}$ on
    pre-Hilbert modules over a fixed $^*$-algebra $\mathcal{D}$.
\end{example}
\begin{example}
    \label{example:ChangeOfBaseRing}
    \emph{(Change of base ring)}
    \\
    Similarly, we can change the base algebra $\mathcal{D}$ in
    $\repC[D](\mathcal{A})$:
    Let $\mathcal{A}$, $\mathcal{D}$ and $\mathcal{D}'$ be
    $^*$-algebras
    and let $\DGDp \in \repC\nolimits_{\mathcal{D}'}(\mathcal{D})$.
    Then $\tensM[D]$ induces a functor
    \begin{equation}
        \label{eq:ChangeOfBaseRing}
        \mathsf{S}_{\mathcal{G}} = \cdot \tensM[D] \DGDp:
        \repC[D](\mathcal{A})
        \longrightarrow \repC\nolimits_{\mathcal{D}'}(\mathcal{A})
    \end{equation}
    defined analogously to \eqref{eq:RieffelInduction}.
\end{example}

A direct consequence of Lemma~\ref{lem:assoc} is that the
following diagram commutes up to natural transformations:
    \begin{equation}
        \label{eq:ReSFCommute}
        \begin{CD}
            \repC[D](\mathcal{A})
            @>\mathsf{S}_{\mathcal{G}}>>
            \repC\nolimits_{\mathcal{D}'}(\mathcal{A}) \\
            @V\mathsf{R}_{\mathcal{E}}VV
            @VV\mathsf{R}_{\mathcal{E}}V \\
            \repC[D](\mathcal{B})
            @>\mathsf{S}_{\mathcal{G}}>>
            \repC\nolimits_{\mathcal{D}'}(\mathcal{B})
        \end{CD}
\end{equation}

\begin{remark}
    \textit{(Rieffel induction for $C^*$-algebras)}
    \label{remark:OriginalRieffel}
    \\
    In the original setting of $C^*$-algebras
    \cite{rieffel:1974a,rieffel:1974b}, Rieffel's construction relates
    categories of $^*$-re\-pre\-sen\-ta\-tions on \textit{Hilbert
      spaces} (in particular, $\mathcal{D}=\ring{C}=\mathbb{C}$), so
    one needs to consider an extra completion of $\BEA \tensM[A] \AHD$
    with respect to the norm induced by
    \eqref{eq:InternalTensorProduct}. Since $^*$-representations of
    $C^*$-algebras on pre-Hilbert spaces are necessarily bounded, this
    completion is canonical, so one recovers Rieffel's original
    construction from this algebraic approach, see
    \cite{bursztyn.waldmann:2001b}. More generally, in this setting,
    $\mathcal{D}$ could be an arbitrary $C^*$-algebra.
\end{remark}

Examples of algebraic Rieffel induction of $^*$-representations in the
setting of formal deformation quantization can be found, e.g., in
\cite{bursztyn.waldmann:2000b,bursztyn.waldmann:2002a}.

\begin{remark}
    \label{remark:ExternalTensor}
    \emph{(External tensor products)}
    \\
    Let $\mathcal{A}_i$ and $\mathcal{B}_i$ be $^*$-algebras over
    $\ring{C}$, $i=1,2$. The tensor products $\mathcal{A} =
    \mathcal{A}_1 \otimes_{\st{\ring{C}}} \mathcal{A}_2$ and
    $\mathcal{B} = \mathcal{B}_1 \otimes_{\st{\ring{C}}}
    \mathcal{B}_2$ are naturally $^*$-algebras.  Let $\mathcal{E}_i$
    be $(\mathcal{B}_i, \mathcal{A}_i)$-bimodules for $i=1,2$ and
    consider the $(\mathcal{B},\mathcal{A})$-bimodule $\mathcal{E} =
    \mathcal{E}_1 \otimes_{\st{\ring{C}}} \mathcal{E}_2$. If each
    $\mathcal{E}_i$ is endowed with an $\mathcal{A}_i$-valued inner
    product $\SP{\cdot,\cdot}_i$, compatible with the
    $\mathcal{B}_i$-action, then we have an inner product
    $\SP{\cdot,\cdot}$ on $\mathcal{E}$, compatible with the
    $\mathcal{B}$-action, uniquely defined by
    \begin{equation}
        \label{eq:ExtTensorSP}
        \SP{x_1 \otimes_{\st{\ring{C}}} x_2, y_1
          \otimes_{\st{\ring{C}}} y_2}
        =
        \SP{x_1, y_1}_1 \otimes_{\st{\ring{C}}} \SP{x_2, y_2}_2
    \end{equation}
    for $x_i, y_i \in \mathcal{E}_i$. We call the inner product
    defined by \eqref{eq:ExtTensorSP} the \textdef{external tensor
      product} of $\SP{\cdot,\cdot}_1$ and $\SP{\cdot,\cdot}_2$.  Just as
      for internal tensor products, if $\SP{\cdot,\cdot}_i$ are
    completely positive, then so is $\SP{\cdot,\cdot}$. The
    construction is also functorial in a sense analogous to
    Lemma~\ref{lemma:MorphInternalTensorMorph}.
\end{remark}

%
%

\section{Strong Morita equivalence}
\label{sec:StrMorEq}

%
%

\subsection{Definition}
\label{sec:SMEDefAndProp}

An $\mathcal{A}$-valued inner product $\SPEA{\cdot,\cdot}$ on an
$\mathcal{A}$-module $\mathcal{E}$ is called \textdef{full} if
\begin{equation}
    \label{eq:Fullness}
    \ring{C} \mbox{-span} \{\SPEA{x,y} \; |\; x,y \in \mathcal{E}\}
    = \mathcal{A}.
\end{equation}
Let $\mathcal{A}$ and $\mathcal{B}$ be $^*$-algebras over $\ring{C}$.
\begin{definition}
    \label{def:strME}
    Let $\BEA$ be a $(\mathcal{B},\mathcal{A})$-bimodule with
    an $\mathcal{A}$-valued inner product $\SPEA{\cdot,\cdot}$ and
    a $\mathcal{B}$-valued inner product $\BSPE{\cdot,\cdot}$. We call
    $(\BEA,\BSPE{\cdot,\cdot},\SPEA{\cdot,\cdot})$ a
    \textdef{$^*$-equivalence bimodule} if
    \begin{enumerate}
    \item $\SPEA{\cdot,\cdot}$ (resp.  $\BSPE{\cdot,\cdot}$) is
        non-degenerate, full and compatible with the
        $\mathcal{B}$-action (resp.  $\mathcal{A}$-action);

    \item For all $x,y,z \in \mathcal{E}$ one has $x\cdot \SPEA{y,z} =
        \BSPE{x,y}\cdot z$;

    \item $\mathcal{B} \cdot \mathcal{E} = \mathcal{E}$ and
        $\mathcal{E} \cdot \mathcal{A}=\mathcal{E}$.
    \end{enumerate}
    If $\SPEA{\cdot,\cdot}$ and $\BSPE{\cdot,\cdot}$ are completely
    positive, then $\BEA$ is called a \textdef{strong equivalence
      bimodule}.
\end{definition}
Whenever the context is clear, we will refer to strong or
$^*$-equivalence bimodules simply as equivalence bimodules.
\begin{definition}
    \label{def:sme}
    Two $^*$-algebras $\mathcal{A}$ and $\mathcal{B}$ are
    \textdef{$^*$-Morita equivalent} (resp. \textdef{strongly Morita
      equivalent}) if there exists a $^*$- (resp. strong)
    $(\mathcal{B},\mathcal{A})$-equivalence bimodule.
\end{definition}

The definition of $^*$-Morita equivalence goes back to Ara
\cite{ara:1999a}.  Since this notion does not involve positivity,
its definition makes sense for ground rings not necessarily of the
form $\ring{C} = \ring{R}(\im)$.
\begin{remark}
    \label{remark:OldFormalME}
    \emph{(Formal Morita equivalence of $^*$-algebras)}
    \\
    In our previous work \cite{bursztyn.waldmann:2001a} we had a more
    technical formulation of strong Morita equivalence for
    $^*$-algebras over $\ring{C}$, called \emph{formal Morita
      equivalence}.  Definition \ref{def:sme}, based on completely
    positive inner products, is conceptually more clear (though, at
    least for unital algebras, it is equivalent to the one of
    \cite{bursztyn.waldmann:2001a}) and yields refinements of the
    results in \cite{bursztyn.waldmann:2001a}.
\end{remark}
\begin{remark}
    \label{remark:GoodOldCstars}
    \emph{(Strong Morita equivalence of $C^*$-algebras)}
    \\
    Rieffel's definition of a strong equivalence bimodule of
    $C^*$-algebras \cite{rieffel:1974b} (see also
    \cite{raeburn.williams:1998a}) is a refinement of
    Definition~\ref{def:strME} involving topological completions which
    do not make sense in a purely algebraic setting. Nevertheless, one
    recovers Rieffel's notion as follows \cite{ara:1999b},
    \cite[Lem.~3.1]{bursztyn.waldmann:2001b}: Two
    $C^*$-algebras are strongly Morita equivalent in Rieffel's sense
    if and only if their minimal dense ideals are strongly Morita
    equivalent (or $^*$-Morita equivalent) in the sense of
    Def.~\ref{def:strME}. In particular, for minimal dense ideals of
    $C^*$-algebras, $^*$- and strong Morita equivalences coincide
    (see Section~\ref{subsec:StrongPicCstar}).
\end{remark}

As we now discuss, $^*$- and strong Morita equivalences are in
fact equivalence relations for a large class of $^*$-algebras. We
start with
\begin{lemma}
    \label{lem:symm}
    The notions of $^*$- and strong Morita equivalences define a
    symmetric relation.
\end{lemma}
For the proof, we just note that if $\BEA$ is a $^*$- (resp.
strong) $(\mathcal{B},\mathcal{A})$-equivalence bimodule, then its
\textit{conjugate bimodule} $\AccEB$ is an $^*$- (resp. strong)
$(\mathcal{A},\mathcal{B})$-equivalence bimodule, see e.g.
\cite[Sect.~5]{bursztyn.waldmann:2001a}.

For reflexivity and transitivity, one needs to be more restrictive.
Recall that an algebra $\mathcal{A}$ is \textdef{non-degenerate} if $a
\in \mathcal{A}$, $\mathcal{A}\cdot a =0$ or $a \cdot \mathcal{A}=0$
implies that $a=0$, and it is \textdef{idempotent} if elements of the
form $a_1 a_2$ span $\mathcal{A}$. The following observation indicates
the importance of these classes of algebras.

Let $\mathcal{A}$ be a $^*$-algebra, and let $\AAA$ be the natural
bimodule induced by left and right multiplications, equipped with the
canonical inner products $\ASP{a,b}= a b^*$ and $\SPA{a,b}=a^*b$.
\begin{lemma}
    \label{lem:reflex}
    The bimodule $\AAA$ is a $^*$- or strong equivalence bimodule if
    and only if $\mathcal{A}$ is non-degenerate and idempotent.
\end{lemma}
The proof is simple: idempotency is equivalent to the canonical
inner products being full and the actions by multiplication being
strongly non-degenerate; non-degeneracy is equivalent to the inner
products being non-degenerate. The inner products are completely
positive by Example~\ref{example:FreeCanonical}.

We therefore restrict ourselves to the class of non-degenerate and
idempotent $^*$-algebras (which contains, in particular, all
unital $^*$-algebras). Within this class, $^*$-Morita equivalence
is transitive \cite{ara:1999a}, hence it is an equivalence
relation. We will show that the same holds for strong Morita
equivalence.

The next result follows from arguments analogous to those in
\cite[Lem.~3.1]{bursztyn.waldmann:2001b}.
\begin{lemma}
    \label{lemma:DegSpacesEqual}
    Let $\mathcal{A}$, $\mathcal{B}$ be non-degenerate and idempotent
    $^*$-algebras, and let $\BEA$ be a bimodule with inner products
    $\SPEA{\cdot,\cdot}$ and $\BSPE{\cdot,\cdot}$ satisfying all the
    properties of Def.~\ref{def:strME}, except for non-degeneracy.
    Then their degeneracy spaces coincide, and the quotient bimodule
    $\mathcal{E} \big/ \mathcal{E}^\bot$, with the induced inner
    products, is a $^*$-equivalence bimodule. If $\SPEA{\cdot,\cdot}$
    and $\BSPE{\cdot,\cdot}$ are completely positive, then the
    quotient bimodule is a strong equivalence bimodule.
\end{lemma}

As a result, within the class of non-degenerate and idempotent
$^*$-algebras, one obtains a refinement of the internal tensor product
$\tensM$ for equivalence bimodules taking into account both inner
products.
\begin{lemma}
    \label{lem:transit}
    Let $\mathcal{A}$, $\mathcal{B}$, $\mathcal{C}$ be non-degenerate
    and idempotent $^*$-algebras and let $\BEA$ and $\CFB$ be $^*$-
    (resp. strong) equivalence bimodules. Then the triple
    \begin{equation}
        \label{eq:EquivBimodTensor}
        \CFB \tensB[B] \BEA :=
        \left(
            \left(\CFB \tensor[B] \BEA\right)
            \big/\left(\CFB \tensor[B] \BEA\right)^\bot,
            \BSPFE{\cdot,\cdot}, \SPFEA{\cdot,\cdot}
        \right)
    \end{equation}
    is a $^*$- (resp. strong) equivalence bimodule.
\end{lemma}
Clearly $\tensB$ satisfies functoriality properties analogous to those
of $\tensM$.
Combining Lemmas
\ref{lem:symm}, \ref{lem:reflex} and \ref{lem:transit}, we obtain:
\begin{theorem}
    \label{theorem:SMEisEquivalence}
    Strong Morita equivalence is an equivalence
    relation within the class of non-degenerate and idempotent $^*$-algebras
    over $\ring{C}$.
\end{theorem}

%
%

\subsection{General properties}
\label{subsec:generalProperties}

Let $\mathcal{A}$ and $\mathcal{B}$ be non-degenerate and idempotent
$^*$-algebras, and let $\Phi: \mathcal{A} \longrightarrow \mathcal{B}$
be a $^*$-isomorphism. A simple check reveals that $\mathcal{B}$, seen
as an $(\mathcal{A},\mathcal{B})$-bimodule via
\begin{equation}
    \label{eq:PhiTwistedBimodule}
    a \cdot_\Phi b \cdot b_1 = \Phi(a)b b_1
\end{equation}
and equipped with the obvious inner products, is a strong
equivalence bimodule. Hence $^*$-isomorphism implies strong Morita
equivalence.

On the other hand, \cite{ara:1999a} shows that $^*$-Morita equivalence
(so also strong Morita equivalence) implies $^*$-isomorphism of
centers. As a result, for commutative (non-degenerate and idempotent)
$^*$-algebras, strong and $^*$-Morita equivalences coincide with the
notion of $^*$-isomorphism.
\begin{remark}
    \label{remark:FiniteRank}
    \emph{(Finite-rank operators)}
    \\
    Let $(\mathcal{E}_{\st{\mathcal{A}}},\SPEA{\cdot,\cdot})$ be an
    inner-product module. The set of ``finite-rank'' operators on
    $\mathcal{E}_{\st{\mathcal{A}}}$, denoted by
    $\mathfrak{F}(\mathcal{E}_{\st{\mathcal{A}}})$, is the
    $\ring{C}$-linear span of operators $\theta_{x,y}$,
    $$
    \theta_{x,y}(z)  := x \cdot \SPEA{y,z},
    $$
    for $x,y,z \in \mathcal{E}$. Note that
    $\theta_{x,y}^*=\theta_{y,x}$ and
    $\mathfrak{F}(\mathcal{E}_{\st{\mathcal{A}}}) \subseteq
    \mathfrak{B}(\mathcal{E}_{\st{\mathcal{A}}})$ is an ideal.

    Within the class of non-degenerate, idempotent $^*$-algebras, an
    alternative description of $^*$-Morita equivalence is given as
    follows \cite{ara:1999a}: if $\mathcal{E}_{\st{\mathcal{A}}}$ is a
    full inner-product module so that
    $\mathcal{E}_{\st{\mathcal{A}}}\cdot\mathcal{A}=\mathcal{E}_{\st{\mathcal{A}}}$,
    then $\FEA$ is a $^*$-equivalence bimodule, with
    $\mathfrak{F}(\mathcal{E}_{\st{\mathcal{A}}})$-valued inner
    product
    \begin{equation}
        \label{eq:thetainner}
        (x,y) \mapsto \theta_{x,y}.
    \end{equation}
    On the other hand, if $\BEA$ is a $^*$-equivalence bimodule, then
    the $\mathcal{B}$-action on $\BEA$ provides a natural
    $^*$-isomorphism
    \begin{equation}
        \label{eq:BisFofEA}
        \mathcal{B}\cong
        \mathfrak{F}(\mathcal{E}_{\st{\mathcal{A}}}).
    \end{equation}
    Under this identification, the $^*$-equivalence bimodule $\BEA$
    becomes $\FEA$.  As a consequence, if
    $\mathcal{E}_{\st{\mathcal{A}}}$ is a pre-Hilbert module with
    $\mathcal{E}_{\st{\mathcal{A}}} \cdot\mathcal{A} =
    \mathcal{E}_{\st{\mathcal{A}}}$ and \eqref{eq:thetainner} is
    completely positive, then $\FEA$ is a strong equivalence bimodule.
\end{remark}

The following is a standard example in Morita theory, see also
\cite[Sect.~6]{bursztyn.waldmann:2001a}.
\begin{example}
    \label{example:MoritaEquivalentAlgebras}
    \emph{(Matrix algebras)}\\~
    Let $\mathcal{A}$ be a non-degenerate and idempotent $^*$-algebra
    over $\ring{C}$. We claim that $\mathcal{A}$ and
    $M_n(\mathcal{A})$ are $^*$- and strongly Morita equivalent.

    First note that $\ring{C}^n$ is a strong
    $(M_n(\ring{C}),\ring{C})$-equivalence bimodule. In fact, since
    $\mathfrak{F}(\ring{C}^n)=M_n(\ring{C})$ and $\ring{C}^n\cdot
    \ring{C} = \ring{C}^n$, by Remark~\ref{remark:FiniteRank} it only
    remains to check that \eqref{eq:thetainner} is completely
    positive. But if $x,y \in \ring{C}^n$, then we can write
    $$
    \theta_{x,y}=\theta_{x,e_1}\theta^*_{y,e_1},
    $$
    where $e_1=(1,0,\ldots,0) \in \ring{C}^n$. So this inner
    product is of the form \eqref{eq:Pinner} (for $m=1$ and
    $P_1(x)=\theta_{x,e_1}$), so it is completely positive.

    By tensoring the equivalence bimodule $\AAA$ with $\ring{C}^n$, it
    follows from Remark~\ref{remark:ExternalTensor} that the canonical
    inner products on $\mathcal{A}^n$ are completely positive. It then
    easily follows that $\mathcal{A}^n$ is a
    $(M_n(\mathcal{A}),\mathcal{A})$-equivalence bimodule.
\end{example}

For unital $^*$-algebras over $\ring{C}$, it follows from the
definitions that strong Morita equivalence implies $^*$-Morita
equivalence, which in turn implies ring-theoretic Morita
equivalence.
In particular, ($^*$- or strong) equivalence bimodules are
finitely generated and projective with respect to both actions.
Using the non-degeneracy of inner products, their compatibility
and fullness, one can verify this property directly by checking
that any $^*$-equivalence bimodule admits a finite hermitian dual
bases. As a consequence, we have

\begin{corollary}
    \label{corollary:unitalcasenondeg}
    If $\mathcal{A}$, $\mathcal{B}$ and $\mathcal{C}$ are unital
    $^*$-algebras and $\CFB$ and $\BEA$ are ($^*$- or strong)
    equivalence bimodules, then the inner product
    \eqref{eq:InternalTensorProduct} on $\CFB \tensor[B] \BEA$ is
    non-degenerate.
\end{corollary}

It follows that the quotient in \eqref{eq:EquivBimodTensor} is
irrelevant. This is always the case for (not necessarily unital)
$C^*$-algebras \cite[Prop.~4.5]{lance:1995a}.

%
%

\subsection{Equivalence of categories of representations}
\label{subsec:EquivalenceCategories}

It is shown in \cite{ara:1999a} that $^*$-Morita equivalence implies
equivalence of categories of (strongly non-degenerate) representations
on inner-product modules. We now recover this result and show that an
analogous statement holds for strong Morita equivalence, generalizing
\cite[Thm.~5.10]{bursztyn.waldmann:2001a}.

The next lemma follows from Lemmas \ref{lem:assoc} and
\ref{lemma:DegSpacesEqual}.
\begin{lemma}
    \label{lemma:naturalisos}
    Let $\mathcal{A}$, $\mathcal{B}$, $\mathcal{C}$ and $\mathcal{D}$
    be non-degenerate and idempotent $^*$-algebras. Let $\CFB$ and
    $\BEA$ be $^*$-equivalence bimodules, and let $(\mathcal{H}, \pi)
    \in \Modrep[D](\mathcal{A})$. Then there are natural isomorphisms
    of inner-product bimodules:
    \begin{align}
        \label{eq:TensBtensMAssoc}
        &\left(\CFB \tensB[B] \BEA\right) \tensM[A] \AHD
        \cong
        \CFB \tensM[B] \left(\BEA \tensM[A] \AHD\right), \\
        &\AAA \tensM[\mathcal{A}] \AHD \cong \AHD \cong
        \AHD \tensM[D] \DDD.
    \end{align}
    As a result, when $\CFB$ and $\BEA$ are strong equivalence
    bimodules, there is a natural equivalence
    \begin{equation}
        \label{eq:RERFREFetc}
        \mathsf{R}_{\mathcal{F}} \circ \mathsf{R}_{\mathcal{E}}
        \cong
        \mathsf{R}_{\mathcal{F} \tensB \mathcal{E}}.
    \end{equation}
\end{lemma}
Using the idempotency and non-degeneracy of $\mathcal{A}$ and
$\mathcal{B}$, one shows
\begin{lemma}
    \label{lemma:InverseBimodules}
    Let $\BEA$ be a ($^*$- or strong) equivalence bimodule.
    If $\AccEB$ is its conjugate bimodule, then
    the following maps are ($^*$- or strong) equivalence bimodule
    isomorphisms:
    \begin{align}
        \label{eq:ccEEisA}
        &\AccEB \tensB[B] \BEA \longrightarrow \AAA,
        \;\;\; \cc{x}\tensB[B] y \mapsto \SPEA{x,y},
        \\\label{eq:ccEEisB}
        &\BEA \tensB[A] \AccEB \longrightarrow \BBB,
        \;\;\; x \tensB[A] \cc{y} \mapsto \BSPE{x,y}.
    \end{align}
\end{lemma}
\begin{corollary}
    \label{corollary:EquiRepTheo}
    Let $\mathcal{A}$, $\mathcal{B}$ and $\mathcal{D}$ be
    non-degenerate and idempotent $^*$-algebras, and let $\BEA$ be a
    strong equivalence bimodule. Then
    \begin{equation}
        \label{eq:RieffelEquivalence}
        \mathsf{R}_{\mathcal{E}}: \RepC[D](\mathcal{A})
        \longrightarrow
        \RepC[D](\mathcal{B})
    \end{equation}
    is an equivalence of categories, with inverse given by
    $\mathsf{R}_{\cc{\mathcal{E}}}$.
\end{corollary}
\begin{remark}
    \label{remark:SFalsoNice}
    Clearly, the functors $\mathsf{S}_{\mathcal{G}}$ satisfy a
    property analogous to \eqref{eq:RERFREFetc}; similarly to
    Corollary \ref{corollary:EquiRepTheo}, an equivalence bimodule
    $\DGDp$ establishes an equivalence of categories
    $\mathsf{S}_{\mathcal{G}}: \RepC[D](\mathcal{A}) \longrightarrow
    \RepC\nolimits_{\mathcal{D}'}(\mathcal{A})$. All these properties
    are direct analogs of the previous constructions by replacing
    tensor products on the left by those on the right.
\end{remark}

Corollary \ref{corollary:EquiRepTheo} recovers the well-knwon theorem
of Rieffel \cite{rieffel:1974b} on the equivalence of categories of
non-degenerate $^*$-representations of strongly Morita equivalent
$C^*$-algebras on Hilbert spaces, see
\cite{bursztyn.waldmann:2001b,ara:1999b}.

%
%

\section{Picard groupoids}
\label{sec:Picard}


In this section, we introduce the Picard groupoids associated with
strong Morita equivalence, in analogy with the groupoid $\Pic$
\cite{benabou:1967a} of invertible bimodules in ring-theoretic
Morita theory \cite{morita:1958a,bass:1968a}.
(See \cite{landsman:2001b,bursztyn.weinstein:2003a:pre} for
related constructions.)

%
%

\subsection{The strong Picard groupoid}
\label{subsec:StrongPicard}

Let $\catstar$ (resp. $\catstarpos$) be the category whose objects
are nondegenerate and idempotent $^*$-algebras over a fixed
$\ring{C}$, morphisms are isomorphism classes of inner-product
(resp.  pre-Hilbert) bimodules and composition is internal tensor
product \eqref{eq:Internal}. (The composition is associative by
Lemma~\ref{lem:assoc}.) We call an inner-product (resp.
pre-Hilbert) bimodule $\BEA$ over $\mathcal{A}$
\textdef{invertible} if its isomorphism class is invertible in
$\catstar$ (resp.  $\catstarpos$). Note that $\BEA$ is invertible
if and only if there exists an inner-product (resp. pre-Hilbert)
bimodule $\AEpB$ over $\mathcal{B}$ together with isomorphisms
\begin{equation}
    \label{eq:invertible}
    \AEpB \tensM[B] \BEA \stackrel{\sim}{\longrightarrow} \AAA,\;\;\;
    \BEA \tensM[A] \AEpB \stackrel{\sim}{\longrightarrow} \BBB.
\end{equation}
\begin{theorem}
    \label{thm:starcat}
    An inner-product (resp. pre-Hilbert) bimodule
    ($\BEA,\SPEA{\cdot,\cdot}$) is invertible if and only if there
    exists a $\mathcal{B}$-valued inner product $\BSPE{\cdot,\cdot}$
    making $(\BEA,\SPEA{\cdot,\cdot},\BSPE{\cdot,\cdot})$ into a
    $^*$-(resp. strong) equivalence bimodule. In particular,
    $^*$-(resp. strong) Morita equivalence coincides with the notion
    of isomorphism in $\catstar$ (resp. $\catstarpos$).
\end{theorem}
This is an algebraic version of a similar result in the framework of
$C^*$-algebras \cite{landsman:2001b,schweizer:1999a:pre}, which we
will recover in Section \ref{subsec:StrongPicCstar}.  We need three
main lemmas to prove the theorem.
\begin{lemma}
    \label{lem:inv}
    Let $(\BEA,\SPEA{\cdot,\cdot})$ be an invertible inner-product
    bimodule. Then $\SPEA{\cdot,\cdot}$ is full and $\mathcal{E}\cdot
    \mathcal{A} = \mathcal{E}$. (By Remark~\ref{remark:FiniteRank},
    $\FEA$ is a $^*$-equivalence bimodule.)
\end{lemma}
\begin{proof}
    Let $\AEpB$ be an inner-product bimodule such that
    \eqref{eq:invertible} holds. The fullness of $\SPEA{\cdot,\cdot}$
    is a simple consequence of the idempotency of $\mathcal{A}$ and
    the first isomorphism of \eqref{eq:invertible}.

    For the second assertion, note that $\BEA \cdot \mathcal{A}
    \subseteq \BEA$ is a $(\mathcal{B},\mathcal{A})$ inner-product
    bimodule. Moreover, using the idempotency of $\mathcal{A}$ and the
    fact that $\mathcal{A}\cdot \AEpB = \AEpB$, it is simple to check
    that $\AEpB$ is still an inverse for $\BEA\cdot \mathcal{A}$.
    By uniqueness of inverses (up to
    isomorphism), we get $\BEA = \BEA \cdot \mathcal{A}$.
\end{proof}
\begin{lemma}
    \label{lem:iso}
    Let $\BEA$ be an invertible inner-product bimodule and let $\FB$
    be an inner-product $\mathcal{B}$-module. Then the natural map
    $$
    \mathsf{S}_{\mathcal{E}}: \mathfrak{B}(\FB) \longrightarrow
    \mathfrak{B}(\FB \tensM[B] \BEA), \;\; T \mapsto T\tensM[B] \id,
    $$
    is an isomorphism.
\end{lemma}
\begin{proof}
    Let $\AEpB$ be as in \eqref{eq:invertible}. Then we have an
    induced map
    $$
    \mathsf{S}_{\mathcal{E}^\prime}:
    \mathfrak{B}(\FB\tensM[B]\BEA) \longrightarrow
    \mathfrak{B}((\FB \tensM[B] \BEA)\tensM[A] \AEpB)\cong
    \mathfrak{B}(\FB),
    $$
    since $ (\FB \tensM[B] \BEA)\tensM[A] \AEpB \cong \FB \tensM[B]
    (\BEA\tensM[A] \AEpB) \cong \FB \tensM[B] \BBB\cong \FB. $ One can
    check that $\mathsf{S}_{\mathcal{E}}$ and
    $\mathsf{S}_{\mathcal{E}^\prime}$ are inverses of each other.
\end{proof}
The next result is an algebraic analog of \cite[Prop.~4.7]{lance:1995a}.
\begin{lemma}
    \label{lem:finite}
    Let $\mathcal{F}_{\st{\mathcal{B}}}$ be an inner-product
    $\mathcal{B}$-module so that $\FB \cdot \mathcal{B}= \mathcal{F}$.
    Let $\mathcal{E}_{\st{\mathcal{A}}}$ be an inner-product
    $\mathcal{A}$-module and $\pi:\mathcal{B}\longrightarrow
    \mathfrak{B}(\mathcal{E}_{\st{\mathcal{A}}})$ be a
    $^*$-homomorphism so that $\pi(\mathcal{B})\subseteq
    \mathfrak{F}(\mathcal{E}_{\st{\mathcal{A}}})$. Then
    $$
    \mathsf{S}_{\mathcal{E}}(\mathfrak{F}(\FB))\subseteq
    \mathfrak{F}(\FB\tensM[B]\BEA).
    $$
\end{lemma}
\begin{proof}
    Suppose $y_1,y_2 \in \FB$ and $b \in \mathcal{B}$. Let $y\tensM[B]
    x \in \FB\tensM[B]\BEA$, and let $\theta_{y_1 \cdot b,y_2} \in
    \mathfrak{F}(\FB)$. Then
    \begin{equation}
        \label{eq:aux}
        \mathsf{S}_{\mathcal{E}}(\theta_{y_1 \cdot b,y_2}) (y\tensM[B] x)
        =
        \theta_{y_1 \cdot b,y_2} y \tensM[B] x
        = y_1 \cdot b\SPFB{y_2,y} \tensM[B] x
        =
        y_1\tensM[B] \pi(b \SPFB{y_2,y}) x.
    \end{equation}
    For each $y \in \FB$, consider the map
    $$
    t_y: \BEA \longrightarrow \FB\tensM[B] \BEA,
    \;\;\; t_y(x)=y\tensM[B] x.
    $$
    Then $t_y \in \mathfrak{B}(\BEA,\FB\tensM[B] \BEA)$, with
    adjoint $t^*_y(y^\prime \tensM[B]
    x^\prime)=\pi(\SPFB{y,y^\prime})x^\prime$. We can rewrite
    \eqref{eq:aux} as
    $$
    \mathsf{S}_{\mathcal{E}}(\theta_{y_1 \cdot b,y_2}) (y\tensM[B]
    x)=t_{y_1}\pi(b)t_{y_2}^*(y\tensM[B] x).
    $$
    Since $\pi(b) \in
    \mathfrak{F}(\mathcal{E}_{\st{\mathcal{A}}})$, it follows that
    $\mathsf{S}_{\mathcal{E}}(\theta_{y_1 \cdot b,y_2}) \in
    \mathfrak{\mathcal{F}(\FB \tensM[B] \BEA})$.


    For a general $\theta_{y_1,y_2}$, we use the condition that $\FB
    \cdot \mathcal{B}=\FB$ to write $y_1 = \sum_{\alpha=1}^k
    y_1^{\alpha} \cdot b_{\alpha}$ and we repeat the argument above.
\end{proof}

We now prove Theorem \ref{thm:starcat} following
\cite{landsman:2001b,schweizer:1999a:pre}.\\
\begin{proof}
    The fact that an equivalence bimodule $\BEA$ is invertible is a
    direct consequence of \eqref{eq:ccEEisA} and
    \eqref{eq:ccEEisB}.

    To prove the other direction, suppose that $\BEA$ is an invertible
    inner-product bimodule, with inverse $\AEpB$.  By Lemma
    \ref{lem:iso}, we have two isomorphisms
    \begin{equation}
        \label{eq:isocomp} \mathfrak{B}(\mathcal{B})
        \stackrel{\mathsf{S}_{\mathcal{E}}}{\longrightarrow}
        \mathfrak{B}(\mathcal{E}_{\st{\mathcal{A}}})
        \stackrel{\mathsf{S}_{\mathcal{E}^\prime}}{\longrightarrow}
        \mathfrak{B}(\mathcal{B}),
    \end{equation}
    whose composition is the identity. Recall that $\mathcal{B} =
    \mathfrak{F}(\mathcal{B}) \subseteq \mathfrak{B}(\mathcal{B})$.
    We claim that
    $$
    \mathsf{S}_{\mathcal{E}^\prime}
    (\mathfrak{F}(\mathcal{E}_{\st{\mathcal{A}}}))
    \subseteq \mathcal{B}.
    $$
    Indeed, for $x_1,x_2 \in \BEA$ and $b \in \mathcal{B}$, we have
    $\mathsf{S}_{\mathcal{E}^\prime}(\theta_{b \cdot x_1, x_2}) =
    b\mathsf{S}_{\mathcal{E}^\prime}(\theta_{x_1, x_2})$, which must
    be in $\mathcal{B}$ since $\mathcal{B}\subset
    \mathfrak{B}(\mathcal{B})$ is an ideal. For a general
    $\theta_{x_1, x_2}$, we use the condition $\mathcal{B}\cdot
    \mathcal{E} = \mathcal{E}$ to write $x_1 =\sum_{\alpha=1}^k
    b_{\alpha} \cdot x_1^{\alpha}$ and apply the same argument. By
    symmetry, it then follows that
    \begin{equation}
        \label{eq:symmetry}
        \mathsf{S}_{\mathcal{E}}
        (\mathfrak{F}(\mathcal{E}^\prime_{\st{\mathcal{B}}}))
        \subseteq \mathcal{A}.
    \end{equation}
    We now claim that
    \begin{equation}
        \label{eq:claim2}
        \mathsf{S}_{\mathcal{E}}(\mathcal{B}) \subseteq
        \mathfrak{F}_{\mathcal{A}}(\mathcal{E}).
    \end{equation}
    By Lemma \ref{lem:inv}, $\FEpB$ is a $^*$-equivalence bimodule.
    Let us consider its conjugate $\BccEpF$. Then
    $$
    \BBB \cong
    \BccEpF \tensMF[\mathfrak{F}(\mathcal{E}^\prime)] \FEpB,
    $$
    and, as a consequence,
    \begin{equation}
        \label{eq:ident1}
        \BEA \cong \BBB \tensM[B] \BEA \cong \BccEpF
        \tensMF[\mathfrak{F}(\mathcal{E}^\prime)]
        (\FEpB \tensM[B] \BEA)
        \cong
        \BccEpF \tensMF[\mathfrak{F}(\mathcal{E}^\prime)] \mathcal{A},
    \end{equation}
    where we regard $\mathcal{A}$ as a left
    $\mathfrak{F}(\mathcal{E}^\prime)$-module via \eqref{eq:symmetry}.

    Since $\BccEpF$ is a $^*$-equivalence bimodule, it follows that
    (see Remark~\ref{remark:FiniteRank}) there is a natural
    identification
    \begin{equation}
        \label{eq:ident2}
        \mathcal{B}
        \cong
        \mathfrak{F}(\overline{\mathcal{E}^\prime}_{\st{\mathfrak{F}(\mathcal{E}^\prime)}}).
    \end{equation}
    Let us now consider the map
    $$
    \mathsf{S}_{\mathcal{A}}: \mathfrak{B}(\BccEpF) \longrightarrow
    \mathfrak{B}(\BccEpF
    \tensMF_{\st{\mathfrak{F}(\mathcal{E^\prime})}}\mathcal{A}).
    $$
    By \eqref{eq:symmetry},
    $\mathfrak{F}(\mathcal{E^\prime}_{\st{\mathcal{B}}})$ acts on
    $\mathcal{A}$ via finite-rank operators; since $\BccEpF\cdot
    \mathfrak{F}(\FEpB) = \BccEpF$, we can apply Lemma
    \ref{lem:finite} and use \eqref{eq:ident1} and \eqref{eq:ident2}
    to conclude that \eqref{eq:claim2} holds.  We can restrict the
    isomorphisms in \eqref{eq:isocomp} to
    $$
    \mathcal{B}
    \stackrel{\mathsf{R}_{\mathcal{E}}}{\longrightarrow}
    \mathfrak{F}(\mathcal{E}_{\st{\mathcal{A}}})
    \stackrel{\mathsf{R}_{\mathcal{E^\prime}}}{\longrightarrow} \mathcal{B},
    $$
    which implies that $\mathsf{R}_{\mathcal{E^\prime}}
    (\mathfrak{F}(\mathcal{E}_{\st{\mathcal{A}}})) = \mathcal{B}$;
    since $\mathsf{R}_{\mathcal{E^\prime}}$ is injective, see
    Lemma~\ref{lem:iso}, there is a natural $^*$-isomorphism
    $\mathcal{B} \cong \mathfrak{F}(\BEA)$, so $\BEA$ is a
    $^*$-equivalence bimodule, again by
    Remark~\ref{remark:FiniteRank}.

    If $\BEA$ and $\AEpB$ are pre-Hilbert bimodules, by uniqueness of
    inverses it follows that
    $$
    \AccEB \cong \AEpB
    $$
    as pre-Hilbert bimodules, so the $\mathcal{B}$-valued inner
    product on $\BEA$ must be completely positive. So $\BEA$ is a
    strong equivalence bimodule.
\end{proof}

The invertible arrows in $\catstar$ (resp. $\catstarpos$) form a
``large'' groupoid $\starPic$ (resp. $\StrPic$), called the
\textdef{$^*$-Picard groupoid} (resp. \textdef{strong Picard
groupoid}). By Theorem \ref{thm:starcat}, orbits of $\starPic$
(resp. $\StrPic$) are $^*$-Morita equivalence (resp. strong Morita
equivalence) classes and isotropy groups are isomorphism classes
of self-$^*$-Morita equivalences (resp. self-strong Morita
equivalences), called \textdef{$^*$-Picard groups} (resp.
\textdef{strong Picard groups}).

If we restrict $\Pic$, $\starPic$ and $\StrPic$ to \textit{unital}
$^*$-algebras over $\ring{C}$, we obtain natural groupoid
homomorphisms
\begin{equation}
    \label{eq:AllThePicMaps}
    \StrPic \longrightarrow \starPic \longrightarrow \Pic,
\end{equation}
covering the identity on the base. On morphisms, the first arrow
``forgets'' the complete positivity of inner products, while the
second just picks the bimodules and ``forgets'' all the extra
structure.

In Section~\ref{sec:strongversusalegbraic}, we will discuss further
conditions on unital $^*$-algebras under which the canonical morphism
\begin{equation}
    \label{eq:ThePicMap}
    \StrPic \longrightarrow \Pic
\end{equation}
is injective and surjective.
\begin{remark}
    \label{remark:WhyNotStarPic}
    The first arrow in \eqref{eq:AllThePicMaps} is generally not
    surjective since a bimodule may have inner products with different
    signatures. For the same reason, the second arrow is not injective
    in general.
\end{remark}

%
%

\subsection{Strong Picard groupoids of $C^*$-algebras}
\label{subsec:StrongPicCstar}

Let $\Cstar$ be the category whose objects are $C^*$-algebras and
morphisms are isomorphism classes of Hilbert bimodules, see e.g.
\cite{landsman:2001b}; the composition is given by
Rieffel's internal tensor product in the $C^*$-algebraic sense. The
groupoid of invertible morphisms in this context will be denoted by
$\StrPic_{\st{C^*}}$. The isotropy groups of $\StrPic_{\st{C^*}}$ are
the Picard groups of $C^*$-algebras as in
\cite{brown.green.rieffel:1977a}.

It is shown in
\cite{schweizer:1999a:pre,landsman:2001b,connes:1994a} that
Rieffel's notion of strong Morita equivalence of $C^*$-algebras
coincides with the notion of isomorphism in $\Cstar$. We will show
how this result can be recovered from Theorem~\ref{thm:starcat}.

For a $C^*$-algebra $\mathcal{A}$, let $\Ped{A}$ be its minimal dense
ideal, also referred to as its \textdef{Pedersen ideal}, see
\cite{pedersen:1979a}. Just as $\mathcal{A}$ itself, $\Ped{A}$ is
non-degenerate and idempotent. If $\mathcal{A}$ and $\mathcal{B}$ are
$C^*$-algebras, let $\cstBEA$ be a Hilbert bimodule (in the
$C^*$-algebraic sense, see e.g. \cite[Def.~3.2]{landsman:2001b}) with
inner product $\SPEcstA{\cdot,\cdot}$, and consider the
$(\Ped{B},\Ped{A})$-bimodule
$$
\mathcal{P}(\cstBEA)
:= \Ped{B} \cdot \cstBEA \cdot \Ped{A}.
$$
\begin{lemma}
    \label{lem:P1}
    The bimodule $\mathcal{P}(\cstBEA)$, together with the
    restriction of $\SPEcstA{\cdot,\cdot}$, is a pre-Hilbert
    $(\Ped{B},\Ped{A})$-bimodule (as in
    Definition~\ref{definition:InnerProductBimodule}).
\end{lemma}
\begin{proof}
    It is clear that $\Ped{B}\cdot \mathcal{P}(\cstBEA) =
    \mathcal{P}(\cstBEA)$, and
    $\SPEcstA{\mathcal{P}(\cstBEA),\mathcal{P}(\cstBEA)} \subseteq
    \Ped{A}$, since $\SPEcstA{\cdot,\cdot}$ is $\mathcal{A}$-linear
    and $\Ped{A}\subseteq \mathcal{A}$ is an ideal.

    For any $n\in \mathbb{N}$, $\mathcal{P}(M_n(\mathcal{A})) =
    M_n(\Ped{A})$. So, by \cite[Lem.~3.2]{bursztyn.waldmann:2001b}, $A
    \in M_n(\Ped{A})^+$ (in the algebraic sense of
    Section~\ref{sec:pos}) if and only if $A \in M_n(\mathcal{A})^+$.
    Since $\SPEcstA{\cdot,\cdot}$ is completely positive, so is its
    restriction to $\mathcal{P}(\cstBEA)$ taking values in $\Ped{A}$.
\end{proof}

The next example shows that $\mathcal{P}(\cstBEA) \subseteq
\cstBEA \cdot \Ped{A}$ is essencial to guarantee that the
restriction of the inner product takes values in $\Ped{A}$.
\begin{example}
    \label{example:restricttoPeds}
    Let $X$ be a locally compact Hausdorff space, and consider
    $\mathcal{A} = C_{\infty}(X)$, the algebra of continuous functions vanishing at $\infty$,
     and $\mathcal{B} =\mathbb{C}$. Then
    $\mathcal{E} = C_{\infty}(X)$ is naturally a
    $(\mathcal{B},\mathcal{A})$-Hilbert bimodule. Since $\mathcal{B}$
    is unital, $\mathcal{B} = \Ped{B}$ and $\Ped{B} \cdot \mathcal{E}
    = \mathcal{E}$. But $\Ped{A} = C_0(X)$ is the algebra of compactly
    supported functions. If $X$ is not compact, then
    $\SPEA{\mathcal{E},\mathcal{E}} \nsubseteq \Ped{A}$.
\end{example}

Let $\mathcal{A},\mathcal{B}$ and $\mathcal{C}$ be $C^*$-algebras.
For Hilbert bimodules $\cstCFB$ and $\cstBEA$, we denote their Rieffel
internal tensor product in the $C^*$-algebraic sense by $\cstCFB
\overline{\otimes}_{\st{\mathcal{B}}} \cstBEA$, see
\cite{rieffel:1974a,rieffel:1974b,raeburn.williams:1998a}. A direct
verification gives
\begin{lemma}
    \label{lem:P2}
    There is a canonical isomorphism $\mathcal{P}(\cstCFB
    \overline{\otimes}_{\st{\mathcal{B}}} \cstBEA) \cong
    \mathcal{P}(\cstCFB) \tensM[B] \mathcal{P}(\cstBEA)$.
\end{lemma}
Let us write $\catstar^{(+)}$ for either $\catstar$ or $\catstarpos$.
With some abuse of notation, it follows from Lemmas \ref{lem:P1} and
\ref{lem:P2} that we can define a functor
\begin{equation}
    \label{eq:funcP}
    \mathcal{P}: \Cstar \longrightarrow \catstar^{(+)}
\end{equation}
as follows: on objects,
$
\mathcal{A} \mapsto \mathcal{P}(\mathcal{A});
$
on morphisms,
$
\mathcal{P}([\cstBEA]) = [\mathcal{P}(\cstBEA)].
$
Here $[\;\;]$ denotes the isomorphism class of a (pre-)Hilbert
bimodule.

Any pre-Hilbert bimodule of the form $\mathcal{E} =
\mathcal{P}(\cstBEA)$ must satisfy $\mathcal{E} \cdot \Ped{A} =
\mathcal{E}$. So the maps on morphisms induced by \eqref{eq:funcP},
$\mathcal{P}: \Mor(\mathcal{A},\mathcal{B}) \longrightarrow
\Mor(\mathcal{P}_{\st{\mathcal{A}}},\mathcal{P}_{\st{\mathcal{B}}})$,
are not surjective in general. However, as we will see, the situation
changes if we restrict $\mathcal{P}$ to morphisms which are
invertible.

Our main analytical tool is the next result, see \cite{ara:1999b}.
\begin{proposition}
    \label{prop:tool}
    Let $\mathcal{A}$ and $\mathcal{B}$ be $C^*$-algebras.
    \begin{itemize}
    \item[i)] If $\PBEPA$ is a strong (or $^*$-)equivalence bimodule
        (as in Def.~\ref{def:strME}), then it can be completed to a
        $C^*$-algebraic strong equivalence bimodule $\cstBEA$ in such
        a way that $\mathcal{P}(\cstBEA) \cong \PBEPA$.

    \item[ii)] If $\cstBFA$ and $\cstBEA$ are $C^*$-algebraic strong
        equivalence bimodules with $\mathcal{P}(\cstBFA) \cong
        \mathcal{P}(\cstBEA)$, then $\cstBFA \cong \cstBEA$.
    \end{itemize}
\end{proposition}
\begin{proof}
    The proof of $i)$ follows from the results in \cite{ara:1999b}.
    Note that $\mathcal{E}_{\st{\Ped{A}}}$ can be completed to a full
    Hilbert $\mathcal{A}$-module
    $\widehat{\mathcal{E}}_{\st{\mathcal{A}}}$, see
    e.g.~\cite[p.~5]{lance:1995a}, so that
    ${_{\st{\mathfrak{K}(\widehat{\mathcal{E}}_{\st{\mathcal{A}}})}}
      \widehat{\mathcal{E}}_{\st{\mathcal{A}}}}$ is a strong
    equivalence bimodule. Here
    $\mathfrak{K}(\widehat{\mathcal{E}}_{\st{\mathcal{A}}})$ denotes
    the ``compact'' operators on
    $\widehat{\mathcal{E}}_{\st{\mathcal{A}}}$. Note that
    $\mathcal{E}$ is naturally an $\mathcal{A}$-module and sits in
    $\widehat{\mathcal{E}}_{\st{\mathcal{A}}}$ as a dense
    $\mathcal{A}$-submodule. So
    $\Ped{B}=\mathfrak{F}(\mathcal{E}_{\st{\mathcal{A}}})$ is dense in
    $\mathfrak{K}(\widehat{\mathcal{E}}_{\st{\mathcal{A}}})$, and
    $\mathcal{B}$ is naturally $^*$-isomorphic to
    $\mathfrak{K}(\widehat{\mathcal{E}}_{\st{\mathcal{A}}})$. So
    $\cstBEA$ is a $C^*$-algebraic strong equivalence bimodule. It
    follows from \cite[Thm.~2.4]{ara:1999b} that any
    $(\Ped{B},\Ped{A})$-$^*$-equivalence bimodule is already a strong
    equivalence bimodule, so the same results hold.

    It follows from \cite{ara:1999b} that
    $$
    \Ped{B}\cdot
    \cstBEA = \cstBEA \cdot \Ped{A} =
    \Ped{B} \cdot \cstBEA \cdot \Ped{A}.
    $$
    Since $\Ped{B} \cdot \PBEPA \cdot \Ped{A} = \PBEPA$, we have
    $\PBEPA \subseteq \Ped{B} \cdot \cstBEA \cdot \Ped{A}$.  On the
    other hand, since $\Ped{B} \subset
    \mathfrak{B}(\widehat{\mathcal{E}}_{\st{\mathcal{A}}})$ is an
    ideal, it follows that $\mathcal{E}\subset \widehat{\mathcal{E}}$
    is
    $\mathfrak{B}(\widehat{\mathcal{E}}_{\st{\mathcal{A}}})$-invariant.
    By \cite[Prop.~1.5]{ara:1999b}, $\Ped{B}\cdot \cstBEA \cdot
    \Ped{A} \subseteq \mathcal{E}$. This implies that
    $\mathcal{P}(\cstBEA) = \PBEPA$.

    Part $ii)$ follows from the fact that $\cstBEA$ is a completion of
    $\mathcal{P}(\cstBEA)$ and any two completions must be isomorphic.
\end{proof}
\begin{corollary}
    \label{corollary:InvHilbertBim}
    A Hilbert bimodule $\cstBEA$ is invertible in $\Cstar$ if and only
    if there exists a $\mathcal{B}$-valued inner product
    $\BSPEcst{\cdot,\cdot}$ so that $(\cstBEA,
    \BSPEcst{\cdot,\cdot},\SPEcstA{\cdot,\cdot})$ is a
    ($C^*$-algebraic) strong equivalence bimodule.  In particular, two
    $C^*$-algebras are strongly Morita equivalent if and only if they
    are isomorphic in $\Cstar$.
\end{corollary}
\begin{proof}
    If $\cstBEA$ is invertible in $\Cstar$, then
    $\mathcal{P}(\cstBEA)$ is invertible in $\catstarpos$.  By
    Theorem~\ref{thm:starcat}, there exists a $\mathcal{B}$-valued
    inner product making $\mathcal{P}(\cstBEA)$ into an equivalence
    bimodule. By part $i)$ of Prop.~\ref{prop:tool}, we can complete
    it to a $C^*$-algebraic strong equivalence bimodule, isomorphic to
    $\cstBEA$ as a Hilbert bimodule.
\end{proof}
\begin{corollary}
    \label{corollary:PicForCstarBij}
    For $C^*$-algebras $\mathcal{A}$ and $\mathcal{B}$,
    \begin{equation}
        \label{eq:PedersenPic}
        \mathcal{P}:
        \StrPic_{\st{C^*}}(\mathcal{B},\mathcal{A})
        \longrightarrow
        \StrPic(\Ped{B},\Ped{A})
    \end{equation}
    is a bijection. As a result, $\StrPic_{\st{C^*}}$ is equivalent to
    $\StrPic$ (or $\starPic$) restricted to Pedersen ideals.
\end{corollary}

It follows that the entire strong Morita theory of $C^*$-algebras is
encoded in the algebraic $\StrPic$. Note that, for unital
$C^*$-algebras, $\mathcal{P}$ is just the identity on objects.

%
%

\section{Strong versus ring-theoretic Picard groupoids}
\label{sec:strongversusalegbraic}

It is shown in \cite{beer:1982a} that unital $C^*$-algebras are
strongly Morita equivalent if and only if they are Morita equivalent
as rings. In \cite{bursztyn.waldmann:2002a}, we have shown that the
same is true for hermitian star products.  In terms of Picard
groupoids, these results mean that $\StrPic$ and $\Pic$, restricted to
unital $C^*$-algebras or to hermitian star products, have the same
orbits.  In this section, we study the morphism $\StrPic
\longrightarrow \Pic$ restricted to unital $^*$-algebras satisfying
additional properties, recovering and refining these results in a
unified way.

%
%

\subsection{A restricted class of unital $^*$-algebras}
\label{subsec:StructBimod}

We consider algebraic conditions which capture some important features
of the functional calculus of $C^*$-algebras. Let $\mathcal{A}$ be a
unital $^*$-algebra over $\ring{C}$. The first property is
\begin{description}
\item[(I)] For all $n\in \mathbb{N}$ and $A \in M_n(\mathcal{A})$,
$\Unit + A^*A$ is invertible.
\end{description}
As a first remark we see that \textbf{(I)} also implies that elements
of the form $\Unit + \sum_{r=1}^k A_r^* A_r$ are invertible in
$M_{n}(\mathcal{A})$, simply by applying \textbf{(I)} to
$M_{nk}(\mathcal{A})$. The relevance of this property is illustrated
by the following result \cite[Thm.~26]{kaplansky:1968a}:
\begin{lemma}
    \label{lemma:Kaplansky}
    Suppose $\mathcal{A}$ satisfies {\rm\textbf{(I)}}. Then any idempotent
    $e = e^2 \in M_n(\mathcal{A})$ is equivalent to a projection $P =
    P^2 = P^* \in M_n(\mathcal{A})$.
\end{lemma}

We also need the following property.
\begin{description}
\item[(II)] For all $n \in \mathbb{N}$, let $P_\alpha \in
    M_n(\mathcal{A})$ be pairwise orthogonal projections, i.e.
    $P_\alpha P_\beta = \delta_{\alpha\beta} P_\alpha$, with $\Unit =
    \sum_\alpha P_\alpha$ and let $H \in M_n(\mathcal{A})^+$ be
    invertible. If $[H, P_\alpha] = 0$, then there exists an
    invertible $U \in M_n(\mathcal{A})$ with $H = U^*U$ and
    $[P_\alpha, U] = 0$.
\end{description}
Most of our results will follow from a condition slightly weaker than
\textbf{(II)}:
\begin{description}
\item[(II$^-$)] For all $n \in \mathbb{N}$, invertible $H \in
    M_n(\mathcal{A})^+$, and projection $P$ with $[P,H] = 0$, there
    exists a $U \in M_n(\mathcal{A})$ with $H = U^*U$ and $[P, U] =
    0$.
\end{description}
On the other hand, our main examples satisfy a stronger version of
\textbf{(II)}:
\begin{description}
\item[(II$^+$)] For all $n \in \mathbb{N}$ and $H \in
    M_n(\mathcal{A})^+$ invertible there exists an invertible $U \in
    M_n(\mathcal{A})$ such that $H = U^*U$, and if $[H, P] = 0$ for a
    projection $P$ then $[U, P] = 0$.
\end{description}

Any unital $C^*$-algebra fulfills \textbf{(I)} and
\textbf{(II$^+$)} by  their functional calculus. In
Section~\ref{subsec:starproducts} we show that the same holds for
hermitian star products. The importance of condition \textbf{(II)}
and its variants lies in the next result.
\begin{lemma}
    \label{lemma:HermModUniqueInnProd}
    Let $\mathcal{A}$ satisfy {\rm\textbf{(II$^-$)}} and let $P = P^2
    = P^* \in M_n(\mathcal{A})$. Then any completely positive and
    strongly non-degenerate $\mathcal{A}$-valued inner product on
    $P\mathcal{A}^n$ is isometric to the canonical one.
\end{lemma}
\begin{proof}
    Given such inner product $\SP{\cdot,\cdot}^\prime$ on
    $P\mathcal{A}^n$, we extend it to the free module $\mathcal{A}^n$
    by taking $(\Unit-P)\mathcal{A}^n$ as orthogonal complement with
    the canonical inner product of $\mathcal{A}^n$ restricted to
    it.  Then the result follows from Lemma~\ref{lem:unique} and the
    fact that the isometry on $\mathcal{A}^n$ commutes with $P$.
\end{proof}

Let $R$ be an arbitrary unital ring.  An idempotent $e=(e_{ij}) \in
M_n(R)$ is called \textdef{full} if the ideal in $R$ generated by
$e_{ij}$ coincides with $R$.  One of the main results of Morita theory
for rings \cite{bass:1968a} is that two unital rings $R$ and $S$ are
Morita equivalent if and only if $S \cong e M_n(R) e$ for some full
idempotent $e$. The next theorem is an analogous result for strong
Morita equivalence.

For a projection $P \in M_n(\mathcal{A})$, we consider
$P\mathcal{A}^n$ equipped with its canonical completely positive inner
product.  Note that $P$ is full if and only if this inner product is
full in the sense of \eqref{eq:Fullness}.
\begin{theorem}
    \label{theorem:StructureOfSMEBimodules}
    Let $\mathcal{A}$, $\mathcal{B}$ be unital $^*$-algebras and let
    $(\BEA, \BSPE{\cdot,\cdot}, \SPEA{\cdot,\cdot})$ be a
    $^*$-equivalence bimodule such that $\SPEA{\cdot,\cdot}$ is
    completely positive. If $\mathcal{A}$ satisfies {\rm\textbf{(I)}}
    and {\rm\textbf{(II$^-$)}} then:
    \begin{enumerate}
    \item There exists a full projection $P = P^2 = P^* \in
        M_n(\mathcal{A})$ such that $\mathcal{E}_{\st{\mathcal{A}}}$
        is isometrically isomorphic to $P\mathcal{A}^n$ as a right
        $\mathcal{A}$-module.
    \item $\mathcal{B}$ is $^*$-isomorphic to $PM_n(\mathcal{A})P$ via
        the left action on $\mathcal{E}_{\st{\mathcal{A}}}$ and the
        $\mathcal{B}$-valued inner product is, under this isomorphism,
        given by the canonical $PM_n(\mathcal{A})P$-valued inner
        product on $P\mathcal{A}^n$.
    \item $\BSPE{\cdot,\cdot}$ is completely positive and hence $\BEA$
        is a strong equivalence bimodule.
    \end{enumerate}
    Conversely, if $P$ is a full projection, then $PM_n(\mathcal{A})P$
    is strongly Morita equivalent to $\mathcal{A}$ via
    $P\mathcal{A}^n$.
\end{theorem}
\begin{proof}
    We know that
    $\mathcal{E}_{\st{\mathcal{A}}}$ is finitely generated and
    projective. By \textbf{(I)} we can find a projection $P$ with
    $\mathcal{E}_{\st{\mathcal{A}}} \cong P \mathcal{A}^n$ and by
    \textbf{(II$^-$)} we can choose the isomorphism to be isometric to
    the canonical inner product, according to
    Lemma~\ref{lemma:HermModUniqueInnProd}, proving the first
    statement.

    Since $\BSPE{\cdot,\cdot}$ is full, the left action map is an
    injective $^*$-homomorphism of $\mathcal{B}$ into
    $\mathfrak{B}(\mathcal{E}_{\st{\mathcal{A}}})$. By compatibility,
    $\BSPE{\cdot,\cdot}$ has to be the canonical one and again by
    fullness we see that $\mathcal{B}$ is $^*$-isomorphic to
    $\mathfrak{B}(\mathcal{E}_{\st{\mathcal{A}}}) \cong P
    M_n(\mathcal{A})P$, proving the second statement.

    Since $\SPEA{\cdot,\cdot}$ is full we find $Px_r, Py_r \in
    P\mathcal{A}^n$ with $\Unit_{\mathcal{A}} = \sum_{r=1}^k
    \SP{Px_r,Py_r}$. Since $\Unit_{\mathcal{A}} =
    \Unit_{\mathcal{A}}^*$ we have
    \[
    \sum\nolimits_r \SP{Px_r + Py_r, Px_r + Py_r}
    =
    \Unit_{\mathcal{A}} + \Unit_{\mathcal{A}} +
    \sum\nolimits_r \SP{Px_r, Px_r} + \sum\nolimits_r \SP{Py_r, Py_r}.
    \]
    By \textbf{(I)} and \textbf{(II$^-$)} we find an invertible $U \in
    \mathcal{A}$ such that for $Pz_r = P(x_r+y_r)U^{-1} \in
    P\mathcal{A}^n$ we have $\Unit_{\mathcal{A}} = \sum_r \SP{Pz_r,
      Pz_r}$.  By compatibility, we get
    \begin{equation}
        \label{eq:BSPECompPos}
        \BSPE{Px,Py} = \sum\nolimits_r \BSPE{Px, Pz_r} \BSPE{Pz_r, Py};
    \end{equation}
    the complete positivity of $\BSPE{\cdot,\cdot}$ now follows from
    Remark~\ref{rem:Pinner}.  This also shows the last statement.
\end{proof}

We use Theorem~\ref{theorem:StructureOfSMEBimodules} to show that
condition \textbf{(I)} and \textbf{(II)} are natural from a
Morita-theoretic point of view.
\begin{proposition}
    \label{proposition:KHstrongHSMEInvariant}
    Conditions {\rm\textbf{(I)}} and {\rm\textbf{(II$^+$)}} (resp.
    {\rm\textbf{(II)}}), together, are strongly Morita invariant.
\end{proposition}
\begin{proof}
    Assume that $\mathcal{A}$ satisfies \textbf{(I)} and
    \textbf{(II$^-$)} and $\mathcal{B}$ is strongly Morita equivalent
    to $\mathcal{A}$.  Then $\mathcal{B} \cong P M_n(\mathcal{A}) P$
    for some full projection $P$. If $B \in \mathcal{B}$, then
    $\Unit_{\mathcal{B}} + B^*B$, viewed as element in $P
    M_n(\mathcal{A})P$, can be extended `block-diagonally' to an
    element of the form
    \begin{equation}
        \label{eq:add}
        \Unit_{M_n(\mathcal{A})} + A^*A
    \end{equation}
    by addition of $\Unit_{M_n(\mathcal{A})} - P$. By \textbf{(I)},
    \eqref{eq:add} has an inverse in $M_n(\mathcal{A})$. By
    \textbf{(II$^-$)}, the inverse is again block-diagonal and hence
    gives an inverse of $\Unit_{\mathcal{B}} + B^*B$. Passing from
    $M_n(\mathcal{A})$ to $M_{nm}(\mathcal{A})$, one obtains the
    invertibility of $\Unit_{M_m(\mathcal{B})} + B^*B$ for $B \in
    M_m(\mathcal{B})$. Hence $\mathcal{B}$ satisfies \textbf{(I)}.

    Assume that $\mathcal{A}$ satisfies \textbf{(II$^+$)},
    and let $H \in \mathcal{B}^+$ be invertible. Then
    $$
    H +
    (\Unit_{M_n(\mathcal{A})} - P) \in M_n(\mathcal{A})^+
    $$
    is still positive and invertible. So there is an invertible $V
    \in M_n(\mathcal{A})$ with $H + (\Unit_{M_n(\mathcal{A})} - P) =
    V^*V$, commuting with $P$ since $H + (\Unit_{M_n(\mathcal{A})} -
    P)$ commutes with $P$. Thus $U = PVP$ satisfies $U^*U = H$.
    Moreover, if $Q \in \mathcal{B}$ is a projection with $[H,Q] = 0$,
    then $PQ = Q = QP$ and hence $Q$ commutes with $H +
    (\Unit_{M_n(\mathcal{A})} - P)$. Thus $V$ commutes with $Q$, and
    hence $U$ commutes with $Q$ as well. For $M_m(\mathcal{B})$, the
    reasoning is analogous.  So $\mathcal{B}$ satisfies
    \textbf{(II$^+$)}.

    An analogous but simpler argument shows the same result for
    \textbf{(II)}.
\end{proof}

%
%

\subsection{From $\StrPic$ to $\Pic$}
\label{subsec:SPictoPic}

Let us consider the groupoid morphism
\begin{equation}
    \label{eq:ThePicMap2}
    \StrPic \longrightarrow \Pic
\end{equation}
from the strong Picard groupoid to the (ring-theoretic) Picard
groupoid.  The next result follows from
Theorem~\ref{theorem:StructureOfSMEBimodules}.
\begin{theorem}
    \label{theorem:PicInjection}
    Within the class of unital $^*$-algebras satisfying
    {\rm\textbf{(I)}} and {\rm\textbf{(II$^-$)}}, the groupoid
    morphism \eqref{eq:ThePicMap2} is injective.
\end{theorem}

For surjectivity, first note that if we define the
\textdef{Hermitian $K_0$-group} of a $^*$-algebra $\mathcal{A}$ as
the Grothendieck group $K^H_0(\mathcal{A})$ of the semi-group of
isomorphism classes of finitely generated projective pre-Hilbert
modules over $\mathcal{A}$ equipped with strongly non-degenerate
inner products, then Lemmas~\ref{lemma:Kaplansky} and
\ref{lemma:HermModUniqueInnProd} imply that if $\mathcal{A}$
satisfies {\rm\textbf{(I)}} and {\rm\textbf{(II$^-$)}}, then the
natural group homomorphism  $K_0^H(\mathcal{A}) \longrightarrow
K_0 (\mathcal{A})$ (forgetting inner products) is an isomorphism.
For Picard groupoids, however, we will see that
\eqref{eq:ThePicMap2} is not generally surjective, even if
\textbf{(I)} and \textbf{(II)} hold. In order to discuss this
surjectivity problem, we consider pairs of $^*$-algebras
satisfying the following rigidity property:

\begin{description}
\item[(III)] Let $\mathcal{A}$ and $\mathcal{B}$ be unital $^*$-algebras,
let $P \in M_n(\mathcal{A})$ be a projection, and consider the $^*$-algebra
$PM_n(\mathcal{A})P$. If $\mathcal{B}$ and $PM_n(\mathcal{A})P$ are isomorphic
as unital algebras, then they are $^*$-isomorphic.
\end{description}

The following are the motivating examples.

\begin{example}
    \label{example:CstarHaveI}
\begin{itemize}

   \item   For unital $C^*$-algebras, condition {\textbf{(III)}} is always satisfied:
           If $\mathcal{A}$ is a $C^*$-algebra, then so is
           $PM_n(\mathcal{A})P$, and {\textbf{(III)}}
           follows from the fact that two $C^*$-algebras which are isomorphic as algebras
           must be $^*$-isomorphic \cite[Thm.~4.1.20]{sakai:1971a}.

   \item   Another class of unital $^*$-algebras satisfying {\textbf{(III)}} is that of
            hermitian star products on a Poisson
           manifold $M$, see Section \ref{sec:applications}. In this case, condition
           {\textbf{(III)}} follows from the more general
           fact that two equivalent star products which are
           compatible with involutions of the form $f \mapsto \overline{f} + o(\lambda)$
           must be $^*$-equivalent, see \cite[Lem.~5]{bursztyn.waldmann:2002a}.
\end{itemize}
\end{example}

For unital algebras $\mathcal{A}$ and $\mathcal{B}$, let us
consider the action of the automorphism group $\Aut(\mathcal{B})$
on the set of morphisms $\Pic(\mathcal{B}, \mathcal{A})$ by
\begin{equation}
    \label{eq:AutActsOnPic}
    (\Phi, [\mathcal{E}]) \mapsto [{\,}_\Phi\mathcal{E}];
\end{equation}
here $\mathcal{E}$ is a $(\mathcal{B},\mathcal{A})$-equivalence
bimodule (in the ring-theoretic sense), $\Phi \in \Aut(\mathcal{B})$
and ${\,}_\Phi\mathcal{E}$ coincides with $\mathcal{E}$ as a
$\ring{C}$-module, but its $(\mathcal{B},\mathcal{A})$-bimodule
structure is given by
$$
b \cdot_\Phi x \cdot a := \Phi(b) \cdot x \cdot a,
$$
see  e.g.~\cite{bass:1968a,bursztyn.waldmann:2002a:pre}.
\begin{proposition}
    \label{proposition:MeetEveryOrbit}
    If $\mathcal{A}$ and $\mathcal{B}$ are unital $^*$-algebras satifying {\rm\textbf{(III)}},
     and if $\mathcal{A}$ satisfies {\rm\textbf{(I)}} and
    {\rm\textbf{(II$^-$)}}, then the composed map
    \begin{equation}
        \label{eq:PicStrongToPic}
        \StrPic(\mathcal{B}, \mathcal{A})
        \longrightarrow
        \Pic(\mathcal{B},\mathcal{A})
        \longrightarrow
        \Pic(\mathcal{B},\mathcal{A}) \big/ \Aut(\mathcal{B})
    \end{equation}
    is onto.
\end{proposition}
\begin{proof}
    Let $\BEA$ be an equivalence bimodule for ordinary Morita
    equivalence.  We know that $\mathcal{E}_{\st{\mathcal{A}}} \cong
    e\mathcal{A}^n$, as right $\mathcal{A}$-modules, for some full
    idempotent $e = e^2 \in M_n(\mathcal{A})$, and $\mathcal{B} \cong
    e M_n(\mathcal{A})e$ as associative algebras via the left action.

    By properties \textbf{(I)} and \textbf{(II$^-$)}, we can replace
    $e$ by a projection $P = P^2 = P^*$ and consider the canonical
    $\mathcal{A}$-valued inner product on $P\mathcal{A}^n$. Then
    $PM_n(\mathcal{A})P$ and $\mathcal{A}$ are strongly Morita
    equivalent via $P\mathcal{A}^n$, see Theorem
    \ref{theorem:StructureOfSMEBimodules}. The identification
    $\mathcal{B} \cong PM_n(\mathcal{A})P$ induces a $^*$-involution
    $^\dag$ on $\mathcal{B}$, possibly different from the original
    one. By assumption, there exists a $^*$-isomorphism
    $$
    \Phi: (\mathcal{B}, {\,}^*)
    \longrightarrow (\mathcal{B}, {\,}^\dag)
    $$
    in such a way that ${\,}_\Phi\mathcal{E}$ becomes a strong
    equivalence bimodule.
\end{proof}
\begin{corollary}
    \label{corollary:RingImpliesStrong}
    Within a class of unital $^*$-algebras satisfying {\rm \textbf{(I)}},
    {\rm \textbf{(II)}} and {\rm \textbf{(III)}}, ring-theoretic Morita equivalence
    implies strong Morita equivalence (so the two notions coincide).
\end{corollary}

Proposition~\ref{proposition:MeetEveryOrbit} is an algebraic
refinement of Beer's result for unital $C^*$-algebras
\cite{beer:1982a}, which is recovered by
Corollary~\ref{corollary:RingImpliesStrong}.

The question to be addressed is when \eqref{eq:ThePicMap2} is
surjective, and not only surjective modulo automorphisms.  The
obstruction for surjectivity is expressed in the next condition.
\begin{description}
\item[(IV)] For any $\Phi \in \Aut(\mathcal{A})$ there is an
    invertible $U \in \mathcal{A}$ such that $\Phi^*\Phi^{-1} =
    \Ad(U^*U)$, where $\Phi^*(a) = \Phi(a^*)^*$.
\end{description}
\begin{lemma}
    \label{lemma:ObstructionSurjectivity}
    Assume that a unital $^*$-algebra $\mathcal{A}$ satisfies {\rm\textbf{(I)}} and
    {\rm\textbf{(II$^-$)}}, and let $\BEA$ be a ring-theoretic
    equivalence bimodule whose class $[\BEA]\in
    \Pic(\mathcal{B},\mathcal{A})$ is in the image of
    \eqref{eq:ThePicMap2}. Then its entire $\Aut(\mathcal{B})$-orbit
    is in the image of \eqref{eq:ThePicMap2} if and only if
    $\mathcal{B}$ satisfies {\rm\textbf{(IV)}}.
\end{lemma}
\begin{proof}
    If the isomorphism class of $\BEA$ is in the image of
    \eqref{eq:ThePicMap2}, then there is a full completely positive
    $\mathcal{A}$-valued inner product $\SPEA{\cdot,\cdot}$ which is
    uniquely determined up to isometry by the right
    $\mathcal{A}$-module structure.

    If $\Phi \in \Aut(\mathcal{B})$, then $[{\,}_\Phi\mathcal{E}]$ is
    in the image of \eqref{eq:ThePicMap2} if and only if there is an
    $\mathcal{A}$-valued inner product $\SP{\cdot,\cdot}^\prime$,
    necessarily isometric to $\SPEA{\cdot,\cdot}$ by
    Lemma~\ref{lemma:HermModUniqueInnProd}, which is compatible with
    the $\mathcal{B}$-action modified by $\Phi$.  In this case, the
    $\mathcal{B}$-valued inner product is determined by compatibility,
    and its complete positivity follows from
    Theorem~\ref{theorem:StructureOfSMEBimodules}. Since there exists
    an invertible $U \in \mathfrak{B}_{\mathcal{A}}(\mathcal{E}) =
    \mathcal{B}$ such that
    $$
    \SP{x,y}^\prime = \SPEA{U \cdot x, U \cdot y},
    $$
    condition \textbf{(IV)} easily follows from the non-degeneracy
    of $\SPEA{\cdot,\cdot}$.
\end{proof}
\begin{corollary}
    \label{corollary:FunnyCondSufNec}
    Let $\mathcal{A}$ and $\mathcal{B}$ be unital $^*$-algebras satisfying
    {\rm{\textbf{(III)}}}, and suppose that $\mathcal{A}$ satisfies properties {\rm\textbf{(I)}}
    and {\rm\textbf{(II$^-$)}}. Then the first map in \eqref{eq:PicStrongToPic}
    is surjective if and only if $\mathcal{B}$ satisfies {\rm\textbf{(IV)}}.
\end{corollary}
\begin{corollary}
    \label{corollary:FunnyConIsSMInv}
    Within a class of unital $^*$-algebras satisfying
    {\rm\textbf{(I)}}, {\rm\textbf{(II$^-$)}}  and {\rm\textbf{(III)}}, property
    {\rm\textbf{(IV)}} is strongly Morita invariant.
\end{corollary}
\begin{example}
    \label{ex:CstarAutos}
    \emph{(The case of $C^*$-algebras)}
    \\
    For a unital $C^*$-algebra $\mathcal{A}$, any automorphism $\Phi
    \in \Aut(\mathcal{A})$ can be uniquely decomposed as
    \begin{equation}
        \label{eq:AutoExpDerStarAuto}
        \Phi = \eu^{\im D} \circ \Psi,
    \end{equation}
    where $\Psi$ is a $^*$-automorphism and $D$ is a $^*$-derivation,
    i.e., a derivation with $D(a^*) = D(a)^*$, see
    \cite[Thm.~7.1]{okayasu:1974a} and
    \cite[Cor.~4.1.21]{sakai:1971a}. In this case, \textbf{(IV)} is
    satisfied if and only if, for any $^*$-derivation $D$, the
    automorphism $\eu^{\im D}$ is inner.

    Let us discuss some concrete examples.  If $\mathcal{A}$ is a
    simple $C^*$-algebra or a $W^*$-algebra, then any automorphism is
    inner, see \cite[Thm.~4.1.19]{sakai:1971a}. So \textbf{(IV)} is
    automatically satisfied, and \eqref{eq:ThePicMap2} is surjective.

    In general, however, there may be automorphisms $\Phi$ with
    $\Phi^* = \Phi^{-1}$ such that $\Phi^2$ is \emph{not} inner, in
    which case \eqref{eq:ThePicMap2} is not surjective.  For example,
    consider the compact operators $\mathfrak{K}(\mathcal{H})$ on a
    Hilbert space $\mathcal{H}$ with countable Hilbert basis $e_n$.
    Define $A = A^* \in \mathfrak{B}(\mathcal{H})$ by
    $$
    A e_{2n} = 2 e_{2n}\;\;\mbox{ and }\;\; A e_{2n+1} = e_{2n+1}.
    $$
    Then $\Ad(A)$ induces an automorphism $\Phi$ of
    $\mathfrak{K}(\mathcal{H}) \oplus \mathbb{C}\Unit$ which satisfies
    $\Phi^* = \Phi^{-1}$ but whose square is not inner: clearly
    $\Ad(A)^2 = \Ad(A^2)$ and there is no $B \in
    \mathfrak{K}(\mathcal{H}) \oplus \mathbb{C}\Unit$ with $\Ad(A^2) =
    \Ad(B^*B)$.
\end{example}

%
%

\section{Hermitian deformation quantization}
\label{sec:applications}

We now show that, just as $C^*$-algebras, hermitian star products can
be treated in the framework of Section
\ref{sec:strongversusalegbraic}. The key observation is that the
properties considered in Section \ref{sec:strongversusalegbraic} are
rigid under deformation quantization.

%
%

\subsection{Hermitian and positive deformations of $^*$-algebras}
\label{subsec:posdef}

Let $\mathcal{A}$ be a $^*$-algebra over $\ring{C}$.  Let
$\boldsymbol{\mathcal{A}} = (\mathcal{A}[[\lambda]], \star)$ be an
associative deformation of $\mathcal{A}$, in the sense of
\cite{gerstenhaber:1964a}.  We call this deformation
\textdef{hermitian} if
$$
(a_1\star a_2)^*=a_2^* \star a_1^*,
$$
for all $a_1,a_2 \in \mathcal{A}$. In this case, $^*$ can be
extended to a $^*$-involution making $\boldsymbol{\mathcal{A}}$ into a
$^*$-algebra over $\ring{C}[[\lambda]]$. Note that
$\ring{C}[[\lambda]]=\ring{R}[[\lambda]](\im)$, and
$\ring{R}[[\lambda]]$ has a natural ordering induced from $\ring{R}$,
see Section~\ref{sec:intro}, so all the notions of positivity of
Section~\ref{sec:pos} make sense for $\boldsymbol{\mathcal{A}}$.
We assume $\lambda$ to be real, so  $\cc{\lambda} = \lambda > 0$.

If $\omega = \sum_{r=0}^\infty \lambda^r \omega_r:
\mathcal{A}[[\lambda]] \longrightarrow \ring{C}[[\lambda]]$ is a
positive $\ring{C}[[\lambda]]$-linear functional with respect to
$\star$, then its classical limit $\omega_0: \mathcal{A}
\longrightarrow \ring{C}$ is a positive $\ring{C}$-linear
functional on $\mathcal{A}.$ We say that a hermitian deformation
$\boldsymbol{\mathcal{A}} = (\mathcal{A}[[\lambda]], \star)$ is
\textdef{positive} \cite[Def.~4.1]{bursztyn.waldmann:2000a} if
every positive linear functional on $\mathcal{A}$ can be
\emph{deformed} into a positive linear functional of
$\boldsymbol{\mathcal{A}}$. In the spirit of complete positivity,
we call a deformation $\boldsymbol{\mathcal{A}}$
\textdef{completely positive} if, for all $n \in \mathbb{N}$, the
$^*$-algebras $M_n(\boldsymbol{\mathcal{A}})$ are positive
deformations of $M_n(\mathcal{A})$. We remark that not all
hermitian deformations are positive \cite{bursztyn.waldmann:2004}.

In the following, we shall consider unital $^*$-algebras and assume
that hermitian deformations preserve the units.

%
%

\subsection{Rigidity of properties \textbf{(I)} and \textbf{(II)}}
\label{subsec:rigid}

The next observation is a direct consequence of the definitions.

\begin{lemma}
    \label{lemma:ComPosDef}
    Let $\boldsymbol{\mathcal{A}}$ be a positive deformation of
    $\mathcal{A}$. If $\boldsymbol{a} = a + o(\lambda) \in
    \boldsymbol{\mathcal{A}}$ is positive, then its classical limit $a
    \in \mathcal{A}$ is also positive.
\end{lemma}

A property of a $^*$-algebra $\mathcal{A}$ is said to be
\textdef{rigid} under a certain type of deformation if any such
deformation satisfies the same property.  Clearly, property
\textbf{(I)} is rigid under hermitian deformations. More
interestingly,
\begin{proposition}
    \label{proposition:Rigidity}
    Property {\rm\textbf{(II)}} is rigid under completely positive
    deformations.
\end{proposition}
\begin{proof}
    Let $\boldsymbol{H} = H + o(\lambda) \in
    M_n(\boldsymbol{\mathcal{A}})^+$ be positive and invertible, and
    let $\boldsymbol{P}_\alpha = P_\alpha + o(\lambda) \in
    M_n(\boldsymbol{\mathcal{A}})$ be pairwise orthogonal projections
    satisfying $\sum_\alpha \boldsymbol{P}_\alpha = \Unit$ and
    $[\boldsymbol{H}, \boldsymbol{P}_\alpha]_\star = 0$.  By
    Lemma~\ref{lemma:ComPosDef}, $H \in M_n(\mathcal{A})$ is positive
    and invertible. Since $[P_\alpha, H] = 0$, by \textbf{(II)} there
    exists an invertible $U \in M_n(\mathcal{A})$ with $H = U^*U$ and
    $[P_\alpha, U] = 0$. In particular, $P_\alpha U P_\alpha \in
    P_\alpha M_n(\mathcal{A}) P_\alpha$ is invertible, with inverse
    $P_\alpha U^{-1} P_\alpha$; here we consider $P_\alpha
    M_n(\mathcal{A}) P_\alpha$ as a unital $^*$-algebra with unit
    $P_\alpha$ as before. Hence
    \begin{equation}
        \label{eq:PalphaHPalpha}
        P_\alpha H P_\alpha = P_\alpha U^* P_\alpha P_\alpha U P_\alpha.
    \end{equation}
    But $\boldsymbol{P}_\alpha \star M_n(\boldsymbol{\mathcal{A}})
    \star \boldsymbol{P}_\alpha$ induces a hermitian deformation
    $\star_\alpha$ of $P_\alpha M_n(\mathcal{A}) P_\alpha$, so we can
    apply \cite[Lem.~2.1]{bursztyn.waldmann:2000b} to deform
    \eqref{eq:PalphaHPalpha}, i.e., there exists an invertible
    $\boldsymbol{U}_\alpha \in \boldsymbol{P}_\alpha \star
    M_n(\boldsymbol{\mathcal{A}}) \star \boldsymbol{P}_\alpha$ such
    that
    \begin{equation}
        \label{eq:DefPalphaHpalpha}
        \boldsymbol{P}_\alpha
        \star \boldsymbol{H}
        \star \boldsymbol{P}_\alpha
        =
        \boldsymbol{P}_\alpha
        \star
        \boldsymbol{U}_\alpha^*
        \star
        \boldsymbol{P}_\alpha
        \star
        \boldsymbol{P}_\alpha
        \star
        \boldsymbol{U}_\alpha
        \star
        \boldsymbol{P}_\alpha.
    \end{equation}
    If we set $\boldsymbol{U} = \sum_\alpha \boldsymbol{U}_\alpha$,
    then it is easy to check that $\boldsymbol{U}$ commutes with the
    $\boldsymbol{P}_\alpha$, is invertible and $\boldsymbol{H} =
    \boldsymbol{U}^* \star \boldsymbol{U}$.
\end{proof}

By only considering the projections $P$ and $(\Unit - P)$, one can
show that property \textbf{(II$^-$)} is rigid under completely
positive deformations as well.

As a consequence,  any completely positive deformation of a $^*$-algebra
$\mathcal{A}$ satisfying \textbf{(I)} and \textbf{(II)} (or
\textbf{(II$^-$)}) also satisfies these properties and those
resulting from them, as discussed in
Section~\ref{sec:strongversusalegbraic}.

%
%

\subsection{Hermitian star products}
\label{subsec:starproducts}

A \textdef{star product} \cite{bayen.et.al:1978a} on a Poisson
manifold $(M,\{\cdot,\cdot\})$ is a formal deformation $\star$ of
$C^\infty(M)$,
$$
f \star g = fg + \sum_{r=1}^\infty \lambda^r C_r(f,g),
$$
for which each $C_r$ is a biddiferential operator and
$$
C_1(f,g)-C_1(g,f)=\im\{f,g\}.
$$
Following Section~\ref{subsec:posdef}, a star product is
\textdef{hermitian} if $\overline{(f\star
  g)}=\overline{g}\star\overline{f}$.

In \cite[Prop.~5.1]{bursztyn.waldmann:2000a}, we proved that any
hermitian star product on a symplectic manifold is a positive
deformation. This turns out to hold much more generally.

\begin{theorem}
    \label{theorem:ComPosStars}
    Any hermitian star product on a Poisson manifold is a completely
    positive deformation.
\end{theorem}
The proof consists of showing that any hermitian star product can
be realized as a subalgebra of a formal Weyl algebra, and then use
the results in the symplectic case \cite{bursztyn.waldmann:2000a},
see \cite{bursztyn.waldmann:2004}.

Since $C^\infty(M)$ satisfies \textbf{(I)} and \textbf{(II)}, we have
\begin{corollary}
    \label{corollary:StarProdKH}
    Hermitian star products on Poisson manifolds satisfy properties
    {\rm\textbf{(I)}} and {\rm\textbf{(II)}}.
\end{corollary}
\begin{corollary}
    \label{corollary:EndE}
    Let $E \longrightarrow M$ be a hermitian vector bundle.  Then any
    hermitian deformation of $\Gamma^\infty(\End(E))$ satisfies
    {\rm\textbf{(I)}} and {\rm\textbf{(II)}}.
\end{corollary}
\begin{proof}
    Any such deformation is strongly Morita equivalent to some
    hermitian star product on $M$, see
    \cite{bursztyn.waldmann:2002a,bursztyn.waldmann:2000b}, so the
    result follows from
    Prop.~\ref{proposition:KHstrongHSMEInvariant}.
\end{proof}

Knowing that hermitian star products are completely positive
deformations, we can use the star exponential to show that they
satisfy a property which is much stronger than \textbf{(II)}.
\begin{proposition}
    \label{proposition:HmuchStronger}
    Let $\star$ be a hermitian star product on $M$. Then  any
    positive invertible $H \in M_n(C^\infty(M)[[\lambda]])^+$ has
    a unique positive invertible $\star$-square root
    $\sqrt[\star]{H}$ such that $[\sqrt[\star]{H}, A]_\star = 0$ if
    and only if $[H, A]_\star = 0$, $A \in
    M_n(C^\infty(M)[[\lambda]])$.  In particular, $\star$ satisfies
    {\rm\textbf{(II$^+$)}}.
\end{proposition}
\begin{proof}
    If $H = H_0 + o(\lambda)$ then, by Lemma~\ref{lemma:ComPosDef},
    $H_0$ is positive in $M_n(C^\infty(M))$ and invertible. This
    implies that $H_0$ has a unique real logarithm $\ln(H_0) \in
    M_n(C^\infty(M))$. Using the star exponential as in
    \cite{bursztyn.waldmann:2002a,bursztyn.waldmann:2002a:pre},
    extended to matrix-valued functions, we conclude that there exists
    a unique real star logarithm
    $
    \Ln(H) = \ln(H_0) + o(\lambda)
    $
     of $H$, whence $\Exp(\Ln(H)) = H$. It follows that $\sqrt[\star]{H} =
    \Exp(\frac{1}{2}\Ln(H))$ has the desired property.
\end{proof}

This shows that many important features of the functional calculus of
$C^*$-algebras are present in formal deformation quantization.

%
%

\subsection{The strong Picard groupoid for star products}
\label{subsec:PicStarProducts}

Since hermitian star products satisfy \textbf{(I)},
\textbf{(II$^+$)} and \textbf{(III)}, it follows that
Thm.~\ref{theorem:PicInjection} and
Prop.~\ref{proposition:MeetEveryOrbit} hold.
\begin{corollary}
    \label{cor:orb}
    For hermitian star products, $\StrPic$ and $\Pic$ have the same
    orbits and the canonical morphism $\StrPic \longrightarrow \Pic$
    is injective.
\end{corollary}

Corollary~\ref{cor:orb} recovers
\cite[Thm.~2]{bursztyn.waldmann:2002a}.  The orbits and isotropy
groups of the Picard groupoid in deformation quantization were
studied in
\cite{bursztyn.waldmann:2002a,bursztyn.waldmann:2002a:pre,jurco.schupp.wess:2002a}.

The next result reveals an interesting similarity between the
structure of the automorphism group of $C^*$-algebras and
hermitian star products, see Example~\ref{ex:CstarAutos}.
\begin{proposition}
    \label{proposition:DecomposeAutos}
    Let $\star$ be a hermitian star product on a Poisson manifold $M$,
    and let $\Phi \in \Aut(\star)$ be an automorphism of $\star$. Then
    there exists a unique $^*$-derivation $D$ and a unique
    $^*$-automorphism $\Psi$ such that
    \begin{equation}
        \label{eq:DecompPhiIntoDandPsi}
        \Phi = \eu^{\im\lambda D} \circ \Psi.
    \end{equation}
\end{proposition}
\begin{proof}
    Writing $\Phi = \sum_{r=0}^\infty \lambda^r\Phi_r$, we know that
    $\Phi_0 = \varphi^*$ is the pull-back by some Poisson
    diffeomorphism $\varphi:M \longrightarrow M$. In particular,
    $\Phi_0 (\cc{f}) = \cc{\Phi_0(f)}$.

    Let us define a new star product $\star'$ by $f \star' g =
    \varphi^*(\varphi_* f \star \varphi_*g)$. Then $\star'$ is
    hermitian and $^*$-isomorphic to $\star$ via $\varphi^*$. If we
    write $\Phi = T \circ \varphi^*$, then $T = \id + o(\lambda)$.
    Hence $\star$ and $\star'$ are equivalent via $T$.

    By \cite[Cor.~4]{bursztyn.waldmann:2002a}, there exists a
    $^*$-equivalence $\tilde{T}$ between $\star$ and $\star'$, so
    $\Psi^{(1)} = \varphi^* \circ \tilde{T}$ is a $^*$-automorphism of
    $\star$ deforming $\varphi^*$. By
    \cite[Lem.~5]{bursztyn.waldmann:2002a}, there is a unique
    derivation $D^{(1)}$ so that $\Phi = \eu^{\im\lambda D^{(1)}}
    \circ \Psi^{(1)}$, and we can write $D^{(1)} = D^{(1)}_1 + \im D^{(1)}_2$, where each
    $D^{(1)}_i$ is a $^*$-derivation.  Now the
    Baker-Campbell-Hausdorff formula defines a derivation $D^{(2)}$ by
    $ \eu^{\im\lambda D} \circ \eu^{\lambda D^{(1)}_2} =
    \eu^{\im\lambda D^{(2)}}, $ in such a way that the imaginary part
    of $D^{(2)}$ is of order $\lambda$. By induction, we can split off
    the $^*$-automorphisms $\eu^{\lambda D_2^{(k)}}$ to obtain
    \eqref{eq:DecompPhiIntoDandPsi}. A simple computation shows the uniqueness of this
    decomposition.
\end{proof}

Using this result, we proceed in total analogy with the case of
$C^*$-algebras. A derivation of a star product $\star$ is
\textdef{quasi-inner} if it is of the form $D = \frac{\im}{\lambda}
\ad(H)$ for some $H \in C^\infty(M)[[\lambda]]$.
\begin{theorem}
    \label{theorem:PicStrongToPicStarProduct}
    Let $\star$, $\star'$ be Morita equivalent hermitian star
    products. Then
    \begin{equation}
        \label{eq:PicStarStar}
        \StrPic(\star, \star') \longrightarrow \Pic(\star, \star')
    \end{equation}
    is a bijection if and only if all derivations of $\star$ are quasi-inner.
\end{theorem}
\begin{proof}
    We know that \eqref{eq:PicStarStar} is injective, and it is
    surjective if and only if any automorphism $\Phi$ of $\star$
    satisfies
    $$
    \cc{\Phi} \circ
    \Phi^{-1} = \Ad(\cc{U} \star U)
    $$
    for some invertible function $U$, see
    Corollary~\ref{corollary:FunnyCondSufNec}.  Using
    \eqref{eq:DecompPhiIntoDandPsi}, this is equivalent to the
    condition that, for any $^*$-derivation $D$, $\eu^{-2\im\lambda D}
    = \Ad(\cc{U} \star U)$.

    Since $\cc{U} \star U = \cc{U_0}U_0 + o(\lambda)$ for some
    invertible $U_0$, we can use the unique real star logarithm
    $\Ln(H)$ of $H = \cc{U} \star U$ to write $$
    \eu^{-2\im\lambda D}
    = \Ad(\Exp(\Ln(H))) = \eu^{\ad(\Ln(H))}.
    $$
    Hence we have the equivalent condition $D = \frac{\im}{\lambda}
    \ad(\frac{1}{2}\Ln(H))$ for any $^*$-derivation. Since any
    derivation can be decomposed into real and imaginary parts, each
    being a $^*$-derivation, the statement follows.
\end{proof}

If $\star$ is a hermitian star product for which Poisson
derivations can be deformed into $\star$-derivations in such a way
that hamiltonian vector fields correspond to quasi-inner
derivations, then $\star$-derivations modulo quasi-inner
derivations are in bijection with formal power series with
coefficients in the first Poisson cohomology, see e.g.
\cite{gutt.rawnsley:1999a,bursztyn.waldmann:2002a:pre}. In this
case, \eqref{eq:PicStarStar} is an isomorphism if and only if the
first Poisson cohomology group vanishes. We recall that any
Poisson manifold admits star products with this property
\cite{cattaneo.felder.tomassini:2002a}, and any symplectic star
product is of this type.
\begin{corollary}
    \label{corollary:PicStrForSomeStar}
    If $\star$ is a hermitian star product on a symplectic manifold
    $M$, then \eqref{eq:PicStarStar} is an isomorphism if and only if
    $\HdR^1(M,\mathbb{C})=\{0\}$.
\end{corollary}

%
%

\begin{footnotesize}

\begin{thebibliography}{10}

\bibitem {ara:1999a}
{\sc Ara, P.: }\newblock {\em {M}orita equivalence for rings with involution}.
\newblock Alg. Rep. Theo.  {\bf 2} (1999), 227--247.

\bibitem {ara:1999b}
{\sc Ara, P.: }\newblock {\em Morita equivalence and Pedersen Ideals}.
\newblock Proc. AMS  {\bf 129}.4 (2000), 1041--1049.

\bibitem {bass:1968a}
{\sc Bass, H.: }\newblock {\em Algebraic ${K}$-theory}.
\newblock W. A. Benjamin, Inc., New York, Amsterdam, 1968.

\bibitem {bayen.et.al:1978a}
{\sc Bayen, F., Flato, M., Fr{{\o}}nsdal, C., Lichnerowicz, A., Sternheimer,
  D.: }\newblock {\em Deformation Theory and Quantization}.
\newblock Ann. Phys.  {\bf 111} (1978), 61--151.

\bibitem {beer:1982a}
{\sc Beer, W.: }\newblock {\em On {M}orita equivalence of nuclear
  ${C}\sp{\ast}$-algebras}.
\newblock J. Pure Appl. Algebra  {\bf 26}.3 (1982), 249--267.

\bibitem {benabou:1967a}
{\sc B{\'e}nabou, J.: }\newblock {\em Introduction to Bicategories}.
\newblock In: {\em Reports of the Midwest Category Seminar},   1--77.
  Springer-Verlag, 1967.

\bibitem {bordemann.waldmann:1998a}
{\sc Bordemann, M., Waldmann, S.: }\newblock {\em Formal GNS Construction and
  States in Deformation Quantization}.
\newblock Commun. Math. Phys.  {\bf 195} (1998), 549--583.

\bibitem {brown.green.rieffel:1977a}
{\sc Brown, L.~G., Green, P., Rieffel, M.: }\newblock {\em Stable Isomorphism
  and Strong Morita Equivalence of {$C^*$}-Algebras}.
\newblock Pacific J. Math.  {\bf 71} (1977), 349--363.

\bibitem {bursztyn.waldmann:2000b}
{\sc Bursztyn, H., Waldmann, S.: }\newblock {\em Deformation Quantization of
  Hermitian Vector Bundles}.
\newblock Lett. Math. Phys.  {\bf 53} (2000), 349--365.

\bibitem {bursztyn.waldmann:2000a}
{\sc Bursztyn, H., Waldmann, S.: }\newblock {\em On Positive Deformations of
  {$^*$}-Algebras}.
\newblock In: {\sc Dito, G., Sternheimer, D. (eds.): }\newblock {\em
  Conf{\'e}rence Mosh{\'e} Flato 1999. Quantization, Deformations, and
  Symmetries}, {\em Mathematical Physics Studies} no. {\bf 22},   69--80.
  Kluwer Academic Publishers, Dordrecht, Boston, London, 2000.

\bibitem {bursztyn.waldmann:2001b}
{\sc Bursztyn, H., Waldmann, S.: }\newblock {\em {$^*$}-Ideals and Formal
  Morita Equivalence of {$^*$}-Algebras}.
\newblock Int. J. Math.  {\bf 12}.5 (2001), 555--577.

\bibitem {bursztyn.waldmann:2001a}
{\sc Bursztyn, H., Waldmann, S.: }\newblock {\em Algebraic Rieffel Induction,
  Formal Morita Equivalence and Applications to Deformation Quantization}.
\newblock J. Geom. Phys.  {\bf 37} (2001), 307--364.

\bibitem {bursztyn.waldmann:2002a:pre}
{\sc Bursztyn, H., Waldmann, S.: }\newblock {\em Bimodule deformations,
  {P}icard groups and contravariant connections}.
\newblock K-theory  {\bf 31} (2004), 1--37.

\bibitem {bursztyn.waldmann:2002a}
{\sc Bursztyn, H., Waldmann, S.: }\newblock {\em The characteristic classes of
  {M}orita equivalent star products on symplectic manifolds}.
\newblock Commun. Math. Phys.  {\bf 228} (2002), 103--121.

\bibitem {bursztyn.waldmann:2004}
{\sc Bursztyn, H., Waldmann, S.: }\newblock {\em Hermitian star
products are completely positive deformations}.
\newblock Preprint \textbf{math.QA/0410350}.

\bibitem {bursztyn.weinstein:2003a:pre}
{\sc Bursztyn, H., Weinstein, A.: }\newblock {\em Picard groups in {P}oisson
  geometry}.
\newblock Moscow Math. J. {\bf 4} (2004), 39--66.


\bibitem {cattaneo.felder.tomassini:2002a}
{\sc Cattaneo, A., Felder, G., Tomassini, L.: }
\newblock {\em From local to global deformation quantization of Poisson manifolds}.
\newblock Duke Math. J.  {\bf 2} (2002), 329--352.

\bibitem {cahen.gutt.dewilde:1980a}
{\sc Cahen, M., Gutt, S., DeWilde, M.: }\newblock {\em Local Cohomology of the
  Algebra of $C^\infty$ Functions on a Connected Manifold}.
\newblock Lett. Math. Phys.  {\bf 4} (1980), 157--167.

\bibitem {connes:1994a}
{\sc Connes, A.: }\newblock {\em Noncommutative Geometry}.
\newblock Academic Press, San Diego, New York, London, 1994.

\bibitem {gerstenhaber:1964a}
{\sc Gerstenhaber, M.: }\newblock {\em On the Deformation of Rings and
  Algebras}.
\newblock Ann. Math.  {\bf 79} (1964), 59--103.

\bibitem {gutt.rawnsley:1999a}
{\sc Gutt, S., Rawnsley, J.: }\newblock {\em Equivalence of star products on a
  symplectic manifold; an introduction to Deligne's {\v{C}}ech cohomology
  classes}.
\newblock J. Geom. Phys.  {\bf 29} (1999), 347--392.

\bibitem {jurco.schupp.wess:2002a}
{\sc Jur{\v{c}}o, B., Schupp, P., Wess, J.: }\newblock {\em Noncommutative Line
  Bundles and Morita Equivalence}.
\newblock Lett. Math. Phys.  {\bf 61} (2002), 171--186.

\bibitem {kaplansky:1968a}
{\sc Kaplansky, I.: }\newblock {\em Rings of operators}.
\newblock W. A. Benjamin, Inc., New York-Amsterdam, 1968.

\bibitem {lance:1995a}
{\sc Lance, E.~C.: }\newblock {\em {H}ilbert {$C^*$}-modules. A toolkit for
  operator algebraists}, vol. 210 in {\em London Mathematical Society Lecture
  Note Series}.
\newblock Cambridge University Press, Cambridge, 1995.

\bibitem {landsman:2001b}
{\sc Landsman, N.~P.: }\newblock {\em Quantized reduction as a tensor product}.
\newblock In: {\sc Landsman, N.~P., Pflaum, M., Schlichenmaier, M. (eds.):
  }\newblock {\em Quantization of Singular Symplectic Quotients},   137--180.
  Birkh{\"a}user, Basel, Boston, Berlin, 2001.

\bibitem {morita:1958a}
{\sc Morita, K.: }\newblock {\em Duality for modules and its applications to
  the theory of rings with minimum condition}.
\newblock Sci. Rep. Tokyo Kyoiku Daigaku Sect. A  {\bf 6} (1958), 83--142.

\bibitem {okayasu:1974a}
{\sc Okayasu, T.: }\newblock {\em Polar Decomposition for Isomorphisms of
  {$C^*$}-Algebras}.
\newblock {T\^{o}hoku} Math. J.  {\bf 26} (1974), 541--554.

\bibitem {pedersen:1979a}
{\sc Pedersen, G.~K.: }\newblock {\em {$C^*$}-Algebras and their Automorphism
  Groups}.
\newblock Academic Press, London, New York, San Francisco, 1979.

\bibitem {raeburn.williams:1998a}
{\sc Raeburn, I., Williams, D.~P.: }\newblock {\em Morita equivalence and
  continuous-trace {$C^*$}-algebras}, vol.~60 in {\em Mathematical Surveys and
  Monographs}.
\newblock American Mathematical Society, Providence, RI, 1998.

\bibitem {rieffel:1974a}
{\sc Rieffel, M.~A.: }\newblock {\em Induced representations of
  {$C^*$}-algebras}.
\newblock Adv. Math.  {\bf 13} (1974), 176--257.

\bibitem {rieffel:1974b}
{\sc Rieffel, M.~A.: }\newblock {\em Morita equivalence for {$C^*$}-algebras
  and {$W^*$}-algebras}.
\newblock J. Pure. Appl. Math.  {\bf 5} (1974), 51--96.

\bibitem {sakai:1971a}
{\sc Sakai, S.: }\newblock {\em {$C^*$}-Algebras and {$W^*$}-Algebras}, vol.~60
  in {\em Ergebnisse der Mathematik und ihrer Grenzgebiete}.
\newblock Springer-Verlag, Berlin, Heidelberg, New York, 1971.

\bibitem{schmuedgen:1990a}
{\sc Schm{\"{u}}dgen, K.:} {\em Unbounded Operator Algebras and Representation
  Theory}, Vol.~37 of {\em Operator Theory: Advances and Applications}.
\newblock Birkh{\"{a}}user, Berlin, 1990.

\bibitem {schweizer:1999a:pre}
{\sc Schweizer, J.: }\newblock {\em Crossed products by equivalence bimodules}.
\newblock Preprint Univ. T{\"u}bingen   (1999).

\bibitem {swan:1962a}
{\sc Swan, R.~G.: }\newblock {\em Vector bundles and projective modules}.
\newblock Trans. Amer. Math. Soc.  {\bf 105} (1962), 264--277.

\bibitem {waldmann:2002a}
{\sc Waldmann, S.: }\newblock {\em On the representation theory of deformation
  quantization}.
\newblock In: {\sc Halbout, G. (eds.): }\newblock {\em Deformation
  quantization}, vol.~1 in {\em IRMA Lectures in Mathematics and Theoretical
  Physics},   107--133. Walter de Gruyter, Berlin, New York, 2002.

\bibitem{waldmann:2003c:pre}
{\sc Waldmann, S.:}\newblock {\em The Picard groupoid in deformation
  quantization}. Preprint \textbf{math.QA/0312118}.
    To appear in Lett. Math. Phys.

\end{thebibliography}

\end{footnotesize}

\end{document}